%
%
%
%
\hsize=5in
\baselineskip=12pt
\vsize=20cm
\parindent=10pt
\pretolerance=40
\predisplaypenalty=0
\displaywidowpenalty=0
\finalhyphendemerits=0
\hfuzz=2pt
\frenchspacing
\footline={\ifnum\pageno=1\else\hfil\tenrm\number\pageno\hfil\fi}
%
%
\input amssym.def
\font\titlef=cmbx12
\font\tenib=cmmib10 
\skewchar\tenib='177
\def\boldfonts{\bf\textfont0=\tenbf}
\font\ninerm=cmr9
\font\ninebf=cmbx9
\font\ninei=cmmi9
\skewchar\ninei='177
\font\ninesy=cmsy9
\skewchar\ninesy='60
\font\nineit=cmti9
\def\reffonts{\baselineskip=0.9\baselineskip
	\textfont0=\ninerm
	\def\rm{\fam0\ninerm}%
	\textfont1=\ninei
  \def\bf{\ninebf}%
	\def\it{\nineit}%
	}
%
%
\def\frontmatter{\vbox{}\vskip1cm\bgroup
	\leftskip=0pt plus1fil\rightskip=0pt plus1fil
	\parindent=0pt
	\parfillskip=0pt
	\pretolerance=10000
	}
\def\endfrontmatter{\egroup\bigskip}
\def\title#1{{\titlef#1\par}}
\def\author#1{\bigskip#1\par}

\def\thanks#1{\footnote{}{\reffonts\rm\noindent#1\hfil}}
\def\section#1\par{\ifdim\lastskip<\bigskipamount\removelastskip\fi
	\penalty-250\bigskip
	\vbox{\leftskip=0pt plus1fil\rightskip=0pt plus1fil
	\parindent=0pt
	\parfillskip=0pt
	\pretolerance=10000{\boldfonts#1}}\nobreak\medskip
	}
\def\emph#1{{\it #1}\/}
\def\proclaim#1. {\medbreak\bgroup{\noindent\bf#1.}\ \it}
\def\endproclaim{\egroup
	\ifdim\lastskip<\medskipamount\removelastskip\medskip\fi}
\newdimen\itemsize
\def\setitemsize#1 {{\setbox0\hbox{#1\ }
	\global\itemsize=\wd0}}
\def\item#1 #2\par{\ifdim\lastskip<\smallskipamount\removelastskip\smallskip\fi
	{\leftskip=\itemsize
	\noindent\hskip-\leftskip
	\hbox to\leftskip{\hfil\rm#1\ }#2\par}\smallskip}
\def\Proof#1. {\ifdim\lastskip<\medskipamount\removelastskip\medskip\fi
	{\noindent\it Proof\if\space#1\space\else\ \fi#1.}\ }
\def\endproof{\hfill\hbox{}\quad\hbox{}\hfill\llap{$\square$}\medskip}
\def\Remark. {\ifdim\lastskip<\medskipamount\removelastskip\medskip\fi
        {\noindent\bf Remark. }}
\def\endremark{\medskip}
%
%
\newcount\citation
\newtoks\citetoks
\def\citedef#1\endcitedef{\citetoks={#1\endcitedef}}
\def\endcitedef#1\endcitedef{}
\def\citenum#1{\citation=0\def\curcite{#1}%
	\expandafter\checkendcite\the\citetoks}
\def\checkendcite#1{\ifx\endcitedef#1?\else
	\expandafter\lookcite\expandafter#1\fi}
\def\lookcite#1 {\advance\citation by1\def\auxcite{#1}%
	\ifx\auxcite\curcite\the\citation\expandafter\endcitedef\else
	\expandafter\checkendcite\fi}
\def\cite#1{\makecite#1,\cite}
\def\makecite#1,#2{[\citenum{#1}\ifx\cite#2]\else\expandafter\clearcite\expandafter#2\fi}
\def\clearcite#1,\cite{, #1]}
%
%
\def\references{\section References\par
	\bgroup
	\parindent=0pt
	\reffonts
	\rm
	\frenchspacing
	\setbox0\hbox{99. }\leftskip=\wd0
	}
\def\endreferences{\egroup}
\newtoks\authtoks
\newif\iffirstauth
\def\checkendauth#1{\ifx\auth#1%
    \iffirstauth\the\authtoks
    \else{} and \the\authtoks\fi,%
  \else\iffirstauth\the\authtoks\firstauthfalse
    \else, \the\authtoks\fi
    \expandafter\nextauth\expandafter#1\fi
	}
\def\nextauth#1,#2;{\authtoks={#1 #2}\checkendauth}
\def\auth#1{\nextauth#1;\auth}
\newif\ifinbook
\newif\ifbookref
\def\nextref#1 {\par\hskip-\leftskip
	\hbox to\leftskip{\hfil\citenum{#1}.\ }%
	\initnextref}
\def\initnextref{\bookreffalse\inbookfalse\firstauthtrue\ignorespaces}
\def\paper#1{{\it#1},}
\def\InBook#1{\inbooktrue in ``#1",}
\def\book#1{\bookreftrue{\it#1},}
\def\journal#1{#1\ifinbook,\fi}
\def\Vol#1{{\bf#1}}
\def\nombre#1{(#1)}
\def\publisher#1{#1,}
\def\Year#1{\ifbookref #1.\else\ifinbook #1,\else(#1)\fi\fi}
\def\Pages#1{\makepages#1.}
\long\def\makepages#1-#2.#3{\ifinbook pp. \fi#1--#2\ifx\par#3.\fi#3}
\def\inRus{{ \rm(in Russian)}}
\def\etransl#1{English translation in \journal{#1}}
%
%
\def\at{\tilde a}
\def\bt{\tilde b}
\def\ct{\tilde c}
\def\xwt{\widetilde x}
\def\thov{\overline\th}
\def\Trhou{T_\rho\mskip1mu u}
\def\Tdu{T_d{\mskip1mu}u}
\newsymbol\square 1003
\newsymbol\bbk 207C 
\newsymbol\curvearrowright 2379
\let\cvar\curvearrowright
\let\ot\otimes
\let\sbs\subset
\newsymbol\smallsetminus 2272
\let\setm\smallsetminus
\let\<\langle
\let\>\rangle
\def\Aut{\mathop{\rm Aut}\nolimits}
\def\GL{GL}
\def\PGL{PGL}
\def\SL{SL}
\def\Id{\mathop{\rm Id}\nolimits}
\def\Im{\mathop{\rm Im}}
\def\Ker{\mathop{\rm Ker}}
\def\Pic{\mathop{\rm Pic}}
\def\Res{\mathop{\rm Res}}
\def\Skl{\mathop{\rm Skl}}
\def\chr{\mathop{\rm char}\nolimits}
\def\op{^{\rm op}}
\def\sgn{\mathop{\rm sgn}\nolimits}
\def\tr{\mathop{\rm tr}\nolimits}
\def\mod#1{\ifinner\mskip8mu(\mathop{\rm mod}#1)
        \else\mskip12mu(\mathop{\rm mod}#1)\fi}
\let\al\alpha
\let\be\beta
\let\ep\varepsilon
\let\ka\kappa
\let\la\lambda
\let\om\omega
\let\ph\varphi
\let\si\sigma
\let\th\theta

\let\La\Lambda
\let\Si\Sigma
\let\Up\Upsilon
\def\frB{{\frak B}}
\def\frS{{\frak S}}
\def\calD{{\cal D}}
\def\calH{{\cal H}}
\def\calI{{\cal I}}
\def\calL{{\cal L}}
\def\bbC{{\Bbb C}}
\def\bbF{{\Bbb F}}
\def\bbP{{\Bbb P}}
\def\bbS{{\Bbb S}}
\def\bbT{{\Bbb T}}
\def\bbZ{{\Bbb Z}}
\let\Rar\Rightarrow
\let\hrar\hookrightarrow

\citedef
Ab-An99
And17
Ar-Sch87
Ar-TV90
Ar-TV91
Bj-Br
Bond-P93
Br-Kn
Ch-WW19
Dub14
Ew-O94
Fi-MS97
Geck-P
Gel-KZ
Gur90
Hai99
Heck-Sch
Lyu87
Man
Ohn99
Ohn05
Pol-P
Ros98
Skr20
Sm96
\endcitedef

\frontmatter

\title{Hecke symmetries: an overview of Frobenius properties}
\author{Serge Skryabin}

\endfrontmatter

\section
Introduction

Hecke symmetries are solutions to the braid version of the quantum Yang-Baxter 
equation satisfying the additional condition $(R-q\cdot\Id)(R+\Id)=0$. 
This notion was introduced by Gurevich \cite{Gur90} as a generalization of 
involutive symmetries considered earlier by Lyubashenko \cite{Lyu87}. Hecke 
symmetries were used to construct nonstandard quantum groups.

With each Hecke symmetry $R$ on a vector space $V$ one associates two graded 
algebras with quadratic defining relations which are viewed as analogs of the 
symmetric and exterior algebras of $V$. We will denote these algebras by 
$\bbS(V,R)$ and $\La(V,R)$, and their homogeneous components by $\bbS_i(V,R)$ 
and $\La_i(V,R)$. A Hecke symmetry $R$ is said to be \emph{even} if $\La(V,R)$ 
is finite dimensional and its highest degree nonzero homogeneous component is 
1-dimensional. In this case Gurevich proved nondegeneracy of bilinear pairings 
between the components of such an algebra, which can be rephrased by saying 
that $\La(V,R)$ is a Frobenius algebra. This fact has important consequences.

However, most of previous work on Hecke symmetries assumed the technical 
condition that the parameter $q$ of the Hecke relation is not a nontrivial 
root of unity. All arguments in \cite{Gur90} exploit projections onto certain 
subspaces of the tensor powers of the base space, while the definition of those 
projections involves division by $q$-integers. Because of that initial setup 
it has gone unnoticed that Gurevich's nondegeneracy result is actually true 
without any restriction on $q$, and the additional condition on Hecke 
symmetries used in \cite{Skr20} is also not needed. We will show a bit more:

\proclaim
Theorem 2.8.
Suppose that\/ $\dim V>1$ and\/ $\dim\La_n(V,R)=1$ for some $n>0$. Then 
$\La_i(V,R)=0$ for all $i>n$ and the algebra\/ $\La(V,R)$ is Frobenius.
\endproclaim

It will be convenient to deal first with an equivalent result for the quadratic 
dual $\bbS(V,R)^!$ of the algebra $\bbS(V,R)$. This will be done in Theorem 2.5. 
In the root of unity case it is not clear whether the algebras $\bbS(V,R)^!$ 
and $\La(V,R)$ always have the same Hilbert series. For this reason we 
distinguish between the algebra $\La(V,R)$ defined as a factor algebra of the 
tensor algebra $\bbT(V)$ and another algebra $\Up(V,R)$ whose homogeneous 
components $\Up^{(k)}$ are subspaces of the tensor powers $V^{\ot k}$. The 
two algebras $\bbS(V,R)^!$ and $\Up(V,R)$ are Frobenius simultaneously.

Every Frobenius algebra admits a so-called Nakayama automorphism determined 
uniquely up to an inner automorphism. In the case of a connected graded 
Frobenius algebra we fix the choice of such an automorphism by requesting it 
to respect the grading. Let $\ph$ be the 1-component of the Nakayama 
automorphism of $\Up(V,R)$, and let $\psi$ be the linear operator on $V$ dual 
to the 1-component of the Nakayama automorphism of $\bbS(V,R)^!$. In section 3 
we describe interrelations between $\ph$, $\psi$ and yet another operator 
$\th\in \GL(V)$ which encodes the braiding
$$
V\ot\Up^{(n)}\to\Up^{(n)}\ot V
$$
arising from $R$. This operator is implicitly present in Gurevich's 
commutation formulas \cite{Gur90, Prop. 5.7}. Here we assume that $\Up^{(n)}$ 
is the last nonzero homogeneous component of $\Up(V,R)$. One of our main 
results is

\proclaim
Theorem 3.8.
The operators $\ph,\psi,\th$ pairwise commute and $\,\psi=q^{-n-1}\ph\,\th^2$. 
Moreover, they lie in the subgroup of $\GL(V)$ consisting of all invertible 
linear operators $\chi:V\to V$ such that $\chi\ot\chi$ commutes with $R$.
\endproclaim

It turns out that the formula $XQ=PY$ in a preprint of Ewen and Ogievetsky 
\cite{Ew-O94, eq. (35)} is the matrix form of the equality $\th\ph=\psi\thov$ 
satisfied with $\thov=q^{n+1}\th^{-1}$. The authors of \cite{Ew-O94} consider 
the 3-dimensional case and call $Q$ and $P$ the characteristic matrices of 
plane and coplane. They correspond to $\ph$ and $\psi$. The matrices $X$ and 
$Y$ are responsible for the commutation law obeyed by the quantum determinant. 
They correspond, up to scalar multiples, to $\th$ and $\thov$. Interrelations 
between those 4 matrices formed a basis for the classification of quantum 
analogs of the group $\GL(3)$.

Actually, Ewen and Ogievetsky postulated that a quantum $\GL(3)$ should have a 
left and a right actions on two graded algebras with the Hilbert series equal 
to that of the polynomial algebra in 3 indeterminates, and the algebra generated 
by the matrix elements of each action should have the same Hilbert series as 
the polynomial algebra in 9 indeterminates. Existence of a Yang-Baxter operator 
satisfying a quadratic relation was deduced from the imposed conditions. 

On the other hand, when $q$ is not a root of unity, it is now known that the 
graded algebras associated with Hecke symmetries have expected Hilbert series 
\cite{Hai99}. Moreover, the graded algebra $\bbS(V,R)$ is Koszul \cite{Gur90}, 
and therefore the Frobenius property of its dual $\bbS(V,R)^!$ ensures that 
$\bbS(V,R)$ is Gorenstein of finite global dimension (see \cite{Bond-P93, 
Prop. 4.1} and \cite{Sm96, Prop. 5.10}). This places $\bbS(V,R)$ in the class 
of algebras introduced by Artin and Schelter \cite{Ar-Sch87} and suggests the 
following

\proclaim
Definition.
A quantum $\GL(3)$ is any Hopf algebra obtained as the Hopf envelope of the 
Faddeev-Reshetikhin-Takhtajan bialgebra associated with a Hecke symmetry $R$ 
on a $3$-dimensional vector space $V$ such that the algebra\/ $\bbS(V,R)$ is 
Artin-Schelter regular of global dimension $3$.
\endproclaim

Two papers of Ohn \cite{Ohn99} and \cite{Ohn05} dealt with the classification 
of quantum analogs of $SL(3)$. The initial postulate there was that the 
representation theory for a quantum $SL(3)$ should be identical to that of the 
ordinary $SL(3)$, including complete reducibility among other things. Each 
quantum $SL(3)$ corresponds to a quantum $\GL(3)$ with central quantum 
determinant. This happens precisely when $\th$ is a scalar operator. All cases 
where $q$ is a nontrivial root of unity are excluded right away \cite{Ohn99, 
Prop. 3.2}. So Ohn's classification does not cover all quantum $\GL(3)$ groups.

The preprint of Ewen and Ogievetsky provides a workable method for determining 
all quantum $\GL(3)$'s and the associated Yang-Baxter operators. It requires 
heavy calculations needed to solve several sets of equations. However, details 
of those calculations are not presented in \cite{Ew-O94}. This makes it 
difficult to verify final conclusions. It seems that the classification of 
quantum groups of $\GL(3)$ type has not been analyzed in later publications.

We would like to know for each Artin-Schelter regular graded quadratic algebra 
$A$ generated by its degree 1 component $A_1=V$ all Hecke symmetries $R$ such 
that $\bbS(V,R)=A$. Twists used in \cite{Ew-O94} are not well suited to answer 
this question as they change the algebra. There are many algebras to consider, 
but this is not the purpose of the present paper. Theorem 5.1 discussed below 
serves as an illustration of general technique.

The Artin-Schelter regular quadratic algebras of global dimension $3$ have 
been subdivided into several types. Those of elliptic type $A$ form a 
particularly interesting class consisting of the Sklyanin algebras 
$\Skl_3(a,b,c)$ for a known set of parameters $a,b,c$. As to the precise 
definition of type $A$ we follow \cite{Ar-TV90}.

\proclaim
Theorem 5.1.
Let $R$ be a Hecke symmetry on a $3$-dimensional vector space $V$ over a field 
of characteristic $\ne2,3$. Then $\bbS(V,R)$ cannot be any Artin-Schelter 
regular algebra of type $A$.
\endproclaim

We will give a full proof of this theorem, allowing arbitrary $q$. In the 
proof we have to examine various possibilities for $\th$. In section 4 we 
describe automorphisms of the algebras $\Skl_3(a,b,c)$ under a suitable 
restriction on $a,b,c$ in terms of the Hessian group $G$ of order 216. The 
auxiliary results there permit us to consider only 4 types of linear operators 
$\th$.

Theorem 5.1 shows that some Artin-Schelter regular quadratic algebras of global 
dimension $3$ are not associated with any quantum $\GL(3)$. One can expect 
abundant examples of such behaviour in higher dimensions. It would be very 
interesting to find a conceptual explanation as to why this happens. The proof 
of Theorem 5.1 is obtained by brute force.

In Manin's approach to quantum groups \cite{Man} one considers the Hopf 
algebra universally coacting on a pair of graded algebras. In \cite{Ch-WW19} 
the authors take two $N$-Koszul Artin-Schelter regular graded algebras 
$A,A'$ and assume a left coaction on $A$ and a right coaction on $A'$. 
Alternatively, one may use right coactions on $A'$ and on the dual of $A$. 
However, without further constraints on the two graded algebras the resulting 
Hopf algebra is likely to be too small to qualify for a true quantum $\GL$. 
Hecke symmetries produce properly compatible pairs of algebras.

In the present paper we consistently use various relations satisfied in the 
Iwahori-Hecke algebras of type $A$. Not much of the theory of Hecke algebras 
is really needed. Several simple facts concerning antisymmetrizers are mentioned 
in section 1, and there we also introduce some notation. A thorough treatment 
of Hecke algebras can be found in \cite{Geck-P}. Combinatorial properties of 
Coxeter groups are discussed in \cite{Bj-Br}. A standard reference on 
quadratic algebras is another book \cite{Pol-P}.

After the first version of this paper had been finished I learned that Abella 
and Andruskiewitsch explained the Frobenius property of $\La(V,R)$ by the fact 
that this algebra has the structure of a braided Hopf algebra \cite{Ab-An99, 
Cor. 3.7.1}. Unlike ordinary Hopf algebras finite-dimensional braided Hopf 
algebras are not always Frobenius. One needs a Hopf module structure on 
the dual space, and its existence can be asserted when the braiding is rigid. 
This condition is usually ensured by considering Hopf algebras in rigid 
monoidal categories, e.g., the category of Yetter-Drinfeld modules for some 
ordinary Hopf algebra. The fact that a braided Hopf algebra is Frobenius was 
stated originally in this setup \cite{Fi-MS97, Cor. 5.8}.

The nondegeneracy result of Gurevich is different in that it actually implies 
that in the case when $\La(V,R)$ is finite dimensional the braiding given by 
the Hecke symmetry $R$ is rigid if and only if the top homogeneous component 
of $\La(V,R)$ has dimension 1. This raises another

\proclaim
Question.
Suppose that $H$ is a finite-dimensional connected graded braided Hopf algebra 
with a $1$-dimensional top homogeneous component. Is then the braiding on $H$ 
necessarily rigid and, as a consequence, $H$ a Frobenius algebra?
\endproclaim

In particular, this question concerns Nichols algebras associated with 
arbitrary braidings (see \cite{And17} or \cite{Heck-Sch}). In our paper we do 
use the quadratic Hecke relation in Lemmas 2.3, 3.1, 3.5, while the main 
results depend on those. I would like to thank N.~Andruskiewitsch for an 
enlightening communication.

\section
1. Antisymmetrizers in the Hecke algebras of type A

Let $\bbk$ be the ground field. The Hecke algebra $\calH_n=\calH_n(q)$ 
of type $A_{n-1}$ with parameter $q\in\bbk$ is presented by generators 
$T_1,\ldots,T_{n-1}$ and relations
$$
\openup1\jot
\displaylines{
T_iT_jT_i=T_jT_iT_j\hbox{\ \ when $|i-j|=1$},\qquad
T_iT_j=T_jT_i\hbox{\ \ when $|i-j|>1$},\cr
(T_i-q)(T_i+1)=0\quad\hbox{for $i=1,\ldots,n-1$}.\cr
}
$$
It has a standard basis $\{T_\si\mid\si\in\frS_n\}$ indexed by elements of the 
symmetric group $\frS_n$. Recall that $\frS_n$ is a Coxeter group with respect 
to its generating set
$$
\frB_n=\{\tau_1,\ldots,\tau_{n-1}\}
$$
of basic transpositions $\tau_i=(i,i+1)$. The \emph{length} $\ell(\si)$ of a 
permutation $\si\in\frS_n$ is the smallest number of factors in the expressions 
of $\si$ as product of basic transpositions. It can be computed as
$$
\ell(\si)=\#\{(i,j)\mid1\le i<j\le n,\ \si(i)>\si(j)\}.
$$
The longest element of $\frS_n$ has length
$$
\ell_n=n(n-1)/2.
$$
So $0\le\ell(\si)\le\ell_n$ for all $\si\in\frS_n$. Denote by $e$ the identity 
permutation. We have $T_e=1$, $\,T_{\tau_i}=T_i$ and
$$
T_\pi T_\si=T_{\pi\si}\quad\hbox{for all $\pi,\si\in\frS_n$ such that 
$\ell(\pi\si)=\ell(\pi)+\ell(\si)$}.
$$

Parabolic subgroups of $\frS_n$ known as \emph{Young subgroups} are generated 
by subsets of $\frB_n$. If $\frS'$ is any Young subgroup, then each left or 
right coset of $\frS'$ in $\frS_n$ contains a unique element of minimal length 
called the \emph{distinguished coset representative}. Distinguished 
representatives of all respective cosets form the subsets
$$
\openup1\jot
\displaylines{
\calD(\frS_n/\frS')=\{\si\in\frS_n\mid\hbox{$\si(i)<\si(i+1)$ for all $i$ 
such that $\tau_i\in\frS'$}\},\cr
\calD(\frS'\backslash\frS_n)=\{\si\in\frS_n\mid\si^{-1}\in\calD(\frS_n/\frS')\} }
$$
of the group $\frS_n$. It is well known that
$$
\openup1\jot
\eqalign{
\ell(d\si)&{}=\ell(d)+\ell(\si)\quad\hbox{for all $\si\in\frS'$ and 
$d\in\calD(\frS_n/\frS')$},\cr
\ell(\si d)&{}=\ell(d)+\ell(\si)\quad\hbox{for all $\si\in\frS'$ and
$d\in\calD(\frS'\backslash\frS_n)$}.
}
$$
For such $\si$ and $d$ we always have, respectively, $T_{d\si}=T_dT_\si$ and 
$T_{\si d}=T_\si T_d\,$.

The algebra $\calH_n$ contains two special elements, the \emph{symmetrizer} 
$\sum_{\si\in\frS_n}T_\si$ and the \emph{antisymmetrizer} 
$$
y_n=\sum_{\si\in\frS_n}(-1)^{\ell(\si)}\,q^{\ell_n-\ell(\si)}\,T_\si\in\calH_n
\qquad\bigl(\ell_n=n(n-1)/2\bigr),
$$
which are interchanged by an involutive automorphism of $\calH_n$. Such 
elements are also defined in the group algebras of Artin's braid groups. Our 
choice in favour of antisymmetrizers is dictated by the convention that 
homogeneous components of the ideal defining $\La(V,R)$ should be 
annihilated by the action of those.

It is known that $y_n$ factorizes as product of partial antisymmetrizers 
in many ways. Let $\frS'$ be any Young subgroup of $\frS_n$. Since each 
element of $\frS_n$ lies in only one left coset and in only one right coset of 
$\frS'$, we get
$$
y_n=y(\frS_n/\frS')\,y(\frS')=y(\frS')\,y(\frS'\backslash\frS_n)
$$
where\quad
$\displaystyle
y(\frS')=\vphantom{\sum}\smash{\sum_{\si\in\frS'}}
(-1)^{\ell(\si)}\,q^{\ell(\frS')-\ell(\si)}\,T_\si,
\qquad
\ell(\frS')=\max\{\ell(\si)\mid\si\in\frS'\},
$
$$
\displaylines{
y(\frS_n/\frS')=\sum_{d\in\calD(\frS_n/\frS')}
(-1)^{\ell(d)}\,q^{\ell_n-\ell(\frS')-\ell(d)}\,T_d\,,\cr
y(\frS'\backslash\frS_n)=\sum_{d\in\calD(\frS'\backslash\frS_n)}
(-1)^{\ell(d)}\,q^{\ell_n-\ell(\frS')-\ell(d)}\,T_d\,.
}
$$
Note that $\ell_n-\ell(\frS')$ is the length of the longest permutations in 
$\calD(\frS_n/\frS')$ and in $\calD(\frS'\backslash\frS_n)$.

Denote by $i\cvar j$ the cycle in $\frS_n$ defined as follows:
$$
\openup1\jot
\eqalign{
i\cvar j&{}=(i,\,i+1,\ldots,j)=\tau_i\tau_{i+1}\cdots\tau_{j-1}\qquad
\hbox{if $i<j$},\cr
i\cvar j&{}=(i,\,i-1,\ldots,j)=\tau_{i-1}\cdots\tau_{j+1}\tau_j\qquad
\hbox{if $i>j$}.
}
$$
In particular, $i\cvar j=j\cvar i=\tau_i$ for $j=i+1$. We set $i\cvar i=e$. In 
this way $i\cvar j\in\frS_n$ are defined for all pairs of integers 
$1\le i,\,j\le n$. The corresponding elements of $\calH_n$ are
$$
T_{i\cvar j}=\cases{T_iT_{i+1}\cdots T_{j-1} & for $i<j$,\cr
\noalign{\smallskip}
1 & for $i=j$,\cr
\noalign{\smallskip}
T_{i-1}\cdots T_{j+1}T_j & for $i>j$.\cr
}
$$
Note that $\ell(i\cvar j)=|j-i|$, and $(i\cvar k)(k\cvar j)=i\cvar j$ whenever 
either $i\le k\le j$ or $i\ge k\ge j$. In this case
$$
T_{i\cvar k}\,T_{k\cvar j}=T_{i\cvar j}\,.
$$

The group $\frS_{n-1}$ is identified with the subgroup of those permutations 
$\si\in\frS_n$ that leave $n$ fixed. Then $\calH_{n-1}$ is identified with the 
subalgebra of $\calH_n$ spanned by all $T_\si$ with $\si\in\frS_{n-1}$. 
Note that
$$
\calD(\frS_n/\frS_{n-1})=\{i\cvar n\mid1\le i\le n\},\qquad
\calD(\frS_{n-1}\backslash\frS_n)=\{n\cvar i\mid1\le i\le n\}.
$$
Taking $\frS'=\frS_{n-1}$ in the previous discussion, we get inductive formulas
$$
y_n=\sum_{i=1}^n(-1)^{n-i}\,q^{i-1}\,T_{i\cvar n}\,y_{n-1}
=y_{n-1}\cdot\sum_{i=1}^n(-1)^{n-i}\,q^{i-1}\,T_{n\cvar i}\,.
$$

Taking $\frS'=\{e,\tau_i\}$ for some $i$, $\,1\le i<n$, we deduce that there 
exist $y',y''\in\calH_n$ such that
$$
y_n=y'\cdot(T_i-q)=(T_i-q)\cdot y''.
$$
Hence $y_n$ is annihilated by $T_i+1$ on both sides, i.e., $T_iy_n=y_nT_i=-y_n$. 
Since this holds for each generator $T_i$ of $\calH_n$, we see that $y_n$ spans 
a 1-dimensional ideal of $\calH_n$, and
$$
T_\si y_n=y_nT_\si=(-1)^{\ell(\si)}y_n\quad\hbox{for all $\si\in\frS_n$}.
$$
We will use standard notation for the $q$-integers and $q$-factorials:
$$
[n]_q=\sum_{i=1}^nq^{i-1},\qquad\quad [n]!_q=\prod_{k=1}^n\,\,[k]_q\,.
$$
There is a combinatorial formula $\sum_{\si\in\frS_n}q^{\ell(\si)}=[n]!_q$ 
proved easily by induction on $n$. It follows that
$$
y_n^2=\sum_{\si\in\frS_n}q^{\ell_n-\ell(\si)}y_n
=\sum_{\si\in\frS_n}q^{\ell(\si)}y_n=[n]!_q\,y_n.
$$
If $[n]!_q\ne0$, then $[n]!_q^{-1}y_n$ is an idempotent in $\calH_n$.

\proclaim
Lemma 1.1.
Given a left $\calH_n$-module $M,$ put 
$$
\Si(M)=\sum_{i=1}^{n-1}\,\{u\in M\mid T_iu=qu\}\quad\hbox{and}\quad
\Up(M)=\bigcap_{i=1}^{n-1}\,(T_i-q)M.
$$
Then $\,y_nM\sbs\Up(M)\,$ and $\,y_n\Si(M)=0$. So the action of $y_n$ on $M$ 
induces a linear map $\,M/\Si(M)\to\Up(M)$.

If\/ $[n]!_q\ne0,$ then this map is bijective, and $\,M=\Up(M)\oplus\Si(M)$.
\endproclaim

\Proof.
The first part of the lemma is a consequence of the factorizations of $y_n$ 
with respect to the 2 element subgroups $\{e,\tau_i\}$ (see above).

Next, the relation $(T_i-q)(T_i+1)=0$ implies that $(T_i+1)u=0$ for all 
$u\in\Up(M)$ and $(T_i+1)u\in\Si(M)$ for all $u\in M$. In other words, 
$\Up(M)$ is a submodule and $M/\Si(M)$ a factor module of $M$ on which each 
$T_i\in\calH_n$ acts as the minus identity operator. Hence $y_nu=[n]!_qu$ for 
all $u\in\Up(M)$ and $y_nu-[n]!_qu\in\Si(M)$ for all $u\in M$. If $[n]!_q\ne0$, we get
$$
\Up(M)=y_nM,\qquad\Si(M)=\{u\in M\mid y_nu=0\},
$$
and the rest is also clear.
\endproof

\proclaim
Lemma 1.2.
Let $M$ be a left $\calH_n$-module, let\/ $\Up(M)$ be as in Lemma 1.1, and let 
$\frS'$ be a Young subgroup of $\frS_n$. Put $\,\calI'=\{j\mid\tau_j\in\frS'\}$. 
Then
$$
y(\frS_n/\frS')\,u\in\Up(M)\quad
\hbox{for each }\,u\in\vphantom{\bigcap}\smash{\bigcap_{j\in\calI'}}(T_j-q)M.
$$
\endproclaim

\Proof.
Put $\calD=\calD(\frS_n/\frS')$. We have to show that the element
$$
y(\frS_n/\frS')\,u=\sum_{d\in\calD}
\,(-1)^{\ell(d)}\,q^{\ell_n-\ell(\frS')-\ell(d)}\,\Tdu
$$
lies in $\,(T_i-q)M\,$ for each $i=1,\ldots,n-1$.

Each permutation $d\in\calD$ has the property that $d(j)<d(j+1)$ for all 
$j\in\calI'$. If $\tau_id\notin\calD$, then $\tau_id(j)>\tau_id(j+1)$ for some 
$j\in\calI'$, but this happens precisely when $i=d(j)$ and $i+1=d(j+1)$. In 
this case $\tau_id=d\tau_j$, and also $T_iT_d=T_{\tau_i d}=T_dT_j$ since 
$\ell(\tau_id)=\ell(d)+1$. Since $\,u=(T_j-q)u'\,$ for some $u'\in M$, we get
$$
\Tdu=(T_i-q)\Tdu'\in(T_i-q)M.
$$
If $\tau_id\in\calD$, then either $\ell(\tau_id)=\ell(d)+1$ or 
$\ell(d)=\ell(\tau_id)+1$. This yields
$$
\openup1\jot
\eqalign{
&y(\frS_n/\frS')\,u=\sum_{\{d\in\calD\mid\tau_id\notin\calD\}}
\,(-1)^{\ell(d)}\,q^{\ell_n-\ell(\frS')-\ell(d)}\,
\Tdu\cr
&\quad{}+\sum_{\{d\in\calD\mid\tau_id\in\calD,\ \ell(d)<\ell(\tau_id)\}}
\,(-1)^{\ell(d)}\,q^{\ell_n-\ell(\frS')-\ell(d)-1}\,
\,(q-T_i)\Tdu\in(T_i-q)M,
}
$$
as required.
\endproof

For two nonnegative integers $k,n$ denote by $\frS_{k,n}$ the Young subgroup 
of $\frS_{k+n}$ consisting of all permutations that leave stable the subset 
$\{1,\ldots,k\}$ and the subset $\{k+1,\ldots,k+n\}$. We assume that $\frS_0$ 
is a trivial group and $\frS_{0,n}=\frS_{n,0}=\frS_n$. Clearly 
$\frS_{k,n}\cong\frS_k\times\frS_n$, and $\frS_{k,n}$ is generated by the 
set $\{\tau_i\in\frB_n\mid i\ne k\}$.

The set $\calD(\frS_{k+n}/\frS_{k,n})$ consists of all \emph{$k,n$-shuffles}. 
These permutations in $\frS_{k+n}$ are increasing on both $\{1,\ldots,k\}$ and 
$\{k+1,\ldots,k+n\}$.

\proclaim
Lemma 1.3.
Let $\rho\in\frS_{k+n}$ be the longest $k,n$-shuffle. Then
$$
\openup2\jot
\displaylines{
\rho=\bigl((1+n)\cvar 1\bigr)\bigl((2+n)\cvar 2\bigr)\cdots
\bigl((k+n)\cvar k\bigr),\cr
T_\rho=T_{(1+n)\cvar 1}\,T_{(2+n)\cvar 2}\,\cdots\,T_{(k+n)\cvar k}\,,\cr
T_iT_\rho=T_{\tau_i\rho}=\cases{T_\rho T_{k+i} & when $1\le i<n,$\cr
T_\rho T_{i-n} & when $n<i<k+n$.
}}
$$
\endproclaim

\Proof.
The permutation $\rho$ maps the set $\{1,\ldots,k\}$ onto $\{1+n,\ldots,k+n\}$ 
and the set $\{k+1,\ldots,k+n\}$ onto $\{1,\ldots,n\}$ preserving the natural 
order on these sets. The displayed expression for $\rho$ is easily checked 
straightforwardly. Also,
$$
\tau_i\rho=\cases{\rho\tau_{k+i} & when $1\le i<n$,\cr
\rho\tau_{i-n} & when $n<i<k+n$.
}
$$
since $\rho$ maps $\{k+i,\,k+i+1\}$ (respectively, $\{i-n,\,i-n+1\}$)
onto $\{i,\,i+1\}$.

The respective equalities in $\calH_{k+n}$ are obtained by observing that 
$\ell(\rho)=kn$ and $\ell(\tau_i\rho)=\ell(\rho)+1$ when $i\ne n$.
\endproof

The partial antisymmetrizers associated with the Young subgroups $\frS_{k,l}$ 
will be used most frequently. We will use abbreviations
$$
y_{k,l}=y(\frS_{k,l}),\qquad
y_{k+l/k,l}=y(\frS_{k+l}/\frS_{k,l}),\qquad
y_{k,l\backslash k+l}=y(\frS_{k,l}\backslash\frS_{k+l}).
$$
Thus for $n=k+l$ we have 
$\displaystyle\ y_n=y_{n/k,l}\cdot y_{k,l}=y_{k,l}\cdot y_{k,l\backslash n}$.

\smallskip
Recall that $\frS_{n-1}$ has been identified with the subgroup $\frS_{n-1,1}$ 
of $\frS_n$. We will write $y_{n/n-1}$ for the element 
$\,y_{n/n-1,1}=y(\frS_n/\frS_{n-1})\in\calH_n$.

\smallskip
In Lemmas 1.4--1.6 we reformulate in our settings several identities known 
for quantum symmetrizers (cf. \cite{Heck-Sch, Cor. 1.8.8}).

\proclaim
Lemma 1.4.
Let $\,1\le k\le n\,$ and $\,c=(n+1)\cvar k\in\frS_{n+1}$. Then
$$
y_{n+1/k,n+1-k}=
q^ky_{n/k,n-k}+(-1)^{n+1-k}y_{n/k-1,n+1-k}\,T_c\,.
$$
\endproclaim

\Proof.
Recall that
$$
y_{n+1/k,n+1-k}=\sum_{d\in\calD(\frS_{n+1}/\frS_{k,n+1-k})}
(-1)^{\ell(d)}\,q^{k(n+1-k)-\ell(d)}\,T_d\,.
$$
For each $d\in\calD(\frS_{n+1}/\frS_{k,n+1-k})$ either $d(n+1)=n+1$ or 
$d(k)=n+1$. In the first of these two cases $d\in\frS_n$, and in fact 
$d\in\calD(\frS_n/\frS_{k,n-k})$. In the second case
$$
d=d'c\quad\hbox{where $d'\in\calD(\frS_n/\frS_{k-1,n+1-k})$},
$$
and we also have $T_d=T_{d'}T_c$ since $c\in\calD(\frS_n\backslash\frS_{n+1})$, 
while $d'\in\frS_n$. Note that $\ell(d)=\ell(d')+n+1-k$. Splitting the above 
sum into two, we arrive at the desired conclusion.
\endproof

For each $h\in\calH_n$ we denote by $h^{(k)}\in\calH_{k+n}$ the image of $h$ 
under the algebra homomorphism $\calH_n\to\calH_{k+n}$ such that $T_i\mapsto 
T_{k+i}$ for $i=1,\ldots,n-1$. We will call $h^{(k)}$ the \emph{$k$ step 
shift} of $h$.

Let similarly $\si^{(k)}\in\frS_{k+n}$ be the image of $\si\in\frS_n$ under the 
group homomorphism $\frS_n\to\frS_{k+n}$ such that $\tau_i\mapsto\tau_{k+i}$ 
for $i=1,\ldots,n-1$. Then $\si^{(k)}$ leaves $1,\ldots,k$ fixed and permutes 
$k+1,\ldots,k+n$ in exactly the same way as $\si$ permutes $1,\ldots,n$. Hence 
$\,\ell(\si^{(k)})=\ell(\si)\,$ and $\,T_\si^{(k)}=T_{\si^{(k)}}\in\calH_{k+n}\,$.

\proclaim
Lemma 1.5.
We have $\displaystyle\ y_{k,l}=y_k{\mskip1mu}y_l^{(k)}$.
\endproclaim

\Proof.
The equality follows from the explicit formula defining $y(\frS_{k,l})$ 
since each permutation in $\frS_{k,l}$ is written uniquely as the product 
$\pi\si^{(k)}$ where $\pi\in\frS_k$ and $\si\in\frS_l$. Clearly 
$\,\ell(\pi\si^{(k)})=\ell(\pi)+\ell(\si)\,$ and 
$\,T_{\pi\si^{(k)}}=T_{\pi}T_{\si^{(k)}}\,$.
\endproof

\proclaim
Lemma 1.6.
We have $\displaystyle\ y_{k+l+m/k+l,m}\cdot y_{k+l/k,l}
=y_{k+l+m/k,l+m}\cdot y_{l+m/l,m}^{(k)}$.
\endproclaim

\Proof.
Both products are equal to $\,y(\frS_{k+l+m}/\frS_{k,l,m})\,$ where 
$\,\frS_{k,l,m}\,$ is the Young subgroup of $\,\frS_{k+l+m}\,$ consisting of 
all permutations that leave stable the 3 subsets
$$
I_1=\{1,\ldots,k\},\quad I_2=\{k+1,\ldots,k+l\},\quad 
I_3=\{k+l+1,\ldots,k+l+m\}.
$$
This is explained as follows. The permutations in 
$\calD(\frS_{k+l+m}/\frS_{k,l,m})$ are increasing on each of these subsets. 
Such a permutation is the product $\pi\si$ where $\pi$ is increasing on 
$I_1\cup I_2$ and on $I_3$, while $\si$ fixes all elements of $I_3$ and is 
increasing on $I_1$ and on $I_2$. Then 
$\pi\in\calD(\frS_{k+l+m}/\frS_{k+l,m})$ and 
$\si\in\calD(\frS_{k+l}/\frS_{k,l})$ where $\frS_{k+l}$ is identified with a 
subgroup of $\frS_{k+l+m}$. We have $\,\ell(\pi\si)=\ell(\pi)+\ell(\si)\,$ and 
$\,T_{\pi\si}=T_\pi T_\si\,$.

On the other hand, each permutation in $\calD(\frS_{k+l+m}/\frS_{k,l,m})$ can 
be written as $\pi\si$ where $\pi$ is increasing on $I_1$ and on $I_2\cup I_3$, 
while $\si$ fixes all elements of $I_1$ and is increasing on $I_2$ and on $I_3$. 
Then $\pi\in\calD(\frS_{k+l+m}/\frS_{k,l+m})$ and $\si=\si'^{\mskip1mu(k)}$ 
for some $\si'\in\calD(\frS_{l+m}/\frS_{l,m})$. Again, 
$\,\ell(\pi\si)=\ell(\pi)+\ell(\si)\,$ and $\,T_{\pi\si}=T_\pi T_\si\,$.
\endproof

\section
2. Nondegeneracy of multiplications

Let $V$ be a finite dimensional vector space over $\bbk$. A \emph{Hecke symmetry} 
on $V$ is a linear operator $R:V\ot V\to V\ot V$ satisfying the braid equation
$$ 
(R\ot\Id_V)(\Id_V\ot\,R)(R\ot\Id_V)=(\Id_V\ot\,R)(R\ot\Id_V)(\Id_V\ot\,R) 
$$
and the quadratic Hecke relation
$$
(R-q\cdot\Id_{V\ot V})(R+\Id_{V\ot V})=0\quad\hbox{where $\,0\ne q\in\bbk$}.
$$
In the case when $q\ne-1$ such an operator is diagonalizable with two 
eigenvalues, but we allow the value $q=-1$ as well.

The Hecke symmetry $R$ gives rise to a representation of $\calH_n$ in the 
$n$th tensor power of $V$ such that $T_i$ acts on $V^{\ot n}$ as the linear 
operator
$$
R_i^{(n)}=\Id_V^{\ot(i-1)}\!\ot\,R\ot\Id_V^{\ot(n-i-1)}\!. 
$$
In this way $V^{\ot n}$ becomes a left $\calH_n$-module. For each $m\ge n$ we 
identify $\calH_n$ with the subalgebra of $\calH_m$ generated by 
$\{T_i\mid i<n\}$, and $V^{\ot m}$ will be always regarded as a left 
$\calH_n$-module with respect to the action
$$
h(a\ot b)=ha\ot b,\qquad h\in\calH_n,\ a\in V^{\ot n},\ b\in V^{\ot(m-n)}.
$$
Note that the $k$ step shift $h^{(k)}\in\calH_{k+n}$ of $h$ (see the 
definition preceding Lemma 1.5) acts as follows:
$$
h^{(k)}(a\ot b)=a\ot hb,\qquad h\in\calH_n,\ a\in V^{\ot k},\ b\in 
V^{\ot m}, m\ge n.
$$

Denote by $\bbT(V)$ the tensor algebra of $V$. It is a graded algebra with 
homogeneous components $\bbT_n(V)=V^{\ot n}$, $\,n\ge0$. The multiplication in 
$\bbT(V)$ is defined by the rule $ab=a\ot b$ for homogeneous elements $a,b$. 
This permits us to omit the sign $\ot$ when writing tensors.

We define the $R$-symmetric algebra $\bbS(V,R)$ and the $R$-skewsymmetric 
algebra $\La(V,R)$ as the factor algebras of $\bbT(V)$ by the graded ideals 
generated, respectively, by the subspaces
$$ 
\Im\,(R-q\cdot\Id)\sbs V^{\ot 2}\quad{\rm and}\quad
\Ker\,(R-q\cdot\Id)\sbs V^{\ot 2}.
$$
Their homogeneous components will be denoted by $\bbS_n(V,R)$ and $\La_n(V,R)$.

Let $I$ be the graded ideal of $\bbT(V)$ such that $\La(V,R)=\bbT(V)/I$. 
Then $I_0=I_1=0$, $\,I_2=\Ker\,(R-q\cdot\Id)$, and
$$
I_n=\sum_{i=1}^{n-1}V^{\ot(i-1)}\!\ot\,I_2\ot V^{\ot(n-i-1)}
=\sum_{i=1}^{n-1}\,\{u\in V^{\ot n}\mid T_iu=qu\}
$$
for $n>2$. Put $\Up^{(0)}=\bbk$, $\Up^{(1)}=V$, $\Up^{(2)}=\Im\,(R-q\cdot\Id)$, 
and
$$
\Up^{(n)}=\bigcap_{i=1}^{n-1}\,(T_i-q)V^{\ot n}
=\bigcap_{i=1}^{n-1}V^{\ot(i-1)}\!\ot\Up^{(2)}\ot V^{\ot(n-i-1)}
$$
for $n>2$. Thus $I_n=\Si(V^{\ot n})$ and $\Up^{(n)}=\Up(V^{\ot n})$ in the 
notation of Lemma 1.1. For each $n$ we identify $\bbT_n(V^*)$ with the dual of 
the vector space $\bbT_n(V)$ using the bilinear pairing
$$
\<f_1\ot\ldots\ot f_n,\,v_1\ot\ldots\ot v_n\>=\prod f_i(v_i).
$$
Under this pairing
$$
(\Up^{(n)})^\perp=\sum_{i=1}^{n-1}(V^*)^{\ot(i-1)}\!\ot(\Up^{(2)})^\perp\ot(V^*)^{\ot(n-i-1)}.
$$
Therefore $\sum_{n=0}^\infty(\Up^{(n)})^\perp$ is a graded ideal of the tensor 
algebra $\bbT(V^*)$ generated by homogeneous elements of degree 2. The factor 
algebra by this ideal is nothing else but the quadratic dual $\bbS(V,R)^!$ of 
$\bbS(V,R)$. Thus the $n$th homogeneous component of $\bbS(V,R)^!$ may be 
identified with the dual space $(\Up^{(n)})^*$.

Let $R^*:(V^*)^{\ot2}\to(V^*)^{\ot2}$ be the linear operator dual to $R$. It is 
a Hecke symmetry on $V^*$. Since $(\Up^{(2)})^\perp=\Ker(R^*-q\cdot\Id)$, we have
$$
\bbS(V,R)^!=\La(V^*,R^*),
$$
and so
$$
\La_n(V^*,R^*)=\bbT_n(V^*)/(\Up^{(n)})^\perp\cong(\Up^{(n)})^*.
$$

We will also consider the subspaces
$$
\Up^{(k,n)}=\Up^{(k)}\ot\Up^{(n)}
=\bigcap_{\{i\,\mid\,1\le i<k+n,\ i\ne k\}}(T_i-q)\,V^{\ot(k+n)}.
$$
If $u\in\Up^{(k,n)}$, then $T_iu=-u$ whenever $1\le i<k+n$ and $i\ne k$. 
Note also that $\Up^{(k+n)}\sbs\Up^{(k,n)}$.

\proclaim
Lemma 2.1.
Assume that $k\ge1$. Let $c=(k+n)\cvar k\in\frS_{k+n},$ and let 
$\rho\in\frS_{k+n}$ be the longest $k,n$-shuffle. Then 
$$
T_c\,\Up^{(k,n)}\sbs\Up^{(k-1,n)}\ot V,\qquad 
T_\rho\,\Up^{(k,n)}\sbs\Up^{(n,k)},\qquad 
y_{k+n/k,n}\Up^{(k,n)}\sbs\Up^{(k+n)}.
$$
\endproclaim

\Proof.
Let $u\in\Up^{(k,n)}$. For $1\le i<k-1$ we have $\tau_ic=c\tau_i$ and 
$\ell(\tau_ic)=\ell(c)+1$, whence $T_iT_c=T_{\tau_ic}=T_cT_i$. Writing 
$u=(T_i-q)w$ for some $w\in V^{\ot(k+n)}$, we get
$$
T_cu=T_c(T_i-q)w=(T_i-q)T_cw.
$$
If $k\le i<k+n-1$, then $\tau_ic=c\tau_{i+1}$ and $\ell(\tau_ic)=\ell(c)+1$, 
which implies that $T_iT_c=T_{\tau_ic}=T_cT_{i+1}$. Since $k<i+1<k+n$, we can 
write $u=(T_{i+1}-q)w$ for some $w\in V^{\ot(k+n)}$, whence
$$
T_cu=T_c(T_{i+1}-q)w=(T_i-q)T_cw.
$$
It follows that\quad
$\displaystyle
T_cu\in\bigcap_{\{i\,\mid\,1\le i<k+n-1,\ i\ne k-1\}}(T_i-q)\,V^{\ot(k+n)}
=\Up^{(k-1,n)}\ot V$.

\medskip
The inclusion $\,T_\rho u\in\Up^{(n,k)}\,$ is proved in a similar way by 
making use of the formulas for $T_iT_\rho$ given in Lemma 1.3.

The inclusion $\,y_{k+n/k,n}\,u\in\Up^{(k+n)}\,$ is a special case of 
Lemma 1.2 applied to the $\calH_{k+n}$-module $M=\bbT_{k+n}(V)$ and the Young 
subgroup $\frS_{k,n}$ of $\frS_{k+n}$.
\endproof

The Hecke symmetry $R$ gives rise to a braiding on a certain monoidal 
subcategory of the category of vector spaces. This was worked out by 
Lyubashenko \cite{Lyu87} in the case of symmetries with parameter $q=1$. A 
special role played by the element $T_\rho\in\calH_{k+n}$ stems from the fact 
that its action on $\bbT_{k+n}(V)$ gives the braiding (commutation in the 
language of \cite{Gur90}) between the spaces $V^{\ot k}$ and $V^{\ot n}$.

In the next lemma we modify the shuffle multiplication of tensors (see Rosso 
\cite{Ros98, Prop. 9}) by using again antisymmetrizers instead of symmetrizers.

\proclaim
Lemma 2.2.
The vector space $\Up(V,R)=\bigoplus_{k=0}^\infty\Up^{(k)}$ is a graded 
associative unital algebra with respect to the multiplication $\star$ defined 
by the formula
$$
a\star b=y_{k+l/k,l}(ab),\qquad a\in\Up^{(k)},\ \,b\in\Up^{(l)}.
$$
The assignments $u\mapsto y_ku$ for $u\in V^{\ot k},$ $k\ge0$ define a 
homomorphism of algebras $\bbT(V)\to\Up(V,R)$ which factors through $\La(V,R)$.  
\endproclaim

\Proof.
By Lemma 2.1 $a\star b\in\Up^{(k+l)}$. Thus the multiplication is well defined, 
and
$$
\eqalign{
(a\star b)\star c&{}=y_{k+l+m/k+l,m}\,y_{k+l/k,l}(abc),\cr
a\star(b\star c)&{}=y_{k+l+m/k,l+m}\,y_{l+m/l,m}^{(k)}(abc)
}
$$
for $a\in\Up^{(k)}$, $b\in\Up^{(l)}$, $c\in\Up^{(m)}$. The associativity law 
$(a\star b)\star c=a\star(b\star c)$ follows from Lemma 1.6. Also, $1\star b=b$ 
and $a\star1=a$ since $y_{k+l/k,l}=1$ whenever either $k=0$ or $l=0$.

By Lemma 1.1 the linear operator by which $y_k$ acts on $V^{\ot k}$ has images 
in $\Up^{(k)}$ and vanishes on the component $I_k$ of the ideal $I\sbs\bbT(V)$ 
which defines $\La(V,R)$. Hence the map $\bbT(V)\to\Up(V,R)$ is well defined 
and factors as required. If $u\in V^{\ot k}$ and $w\in V^{\ot l}$, then 
$$
(y_ku)\star(y_lw)=y_{k+l/k,l}\,y_ky_l^{(k)}(uw)=y_{k+l/k,l}\,y_{k,l}(uw)=y_{k+l}(uw)
$$
by Lemma 1.5.
\endproof

\proclaim
Lemma 2.3.
Suppose that $\Up^{(n+1)}=0$ for some $n>0$. Let $1\le k\le n,$ and let 
$\rho\in\frS_{k+n}$ be the longest $k,n$-shuffle. Then
$$
y_{n/k,n-k}\,u=(-1)^{kn-k}\,q^{-k(k+1)/2}\,\Trhou
\quad\hbox{for all $\,u\in\Up^{(k,n)}$}.
$$
In the case $k=n$ this yields $\,\Trhou=q^{n(n+1)/2}\,u\,$ for $u\in\Up^{(n,n)}$.
\endproclaim

\Proof.
Let $u\in\Up^{(k,n)}$. Note that $\Up^{(k,n)}\sbs\Up^{(k,n+1-k)}\ot\bbT_{k-1}(V)$. 
Therefore 
$$
y_{n+1/k,n+1-k}\,u\in\Up^{(n+1)}\ot\bbT_{k-1}(V)=0
$$
by Lemma 2.1. In view of Lemma 1.4 this can be rewritten as
$$
y_{n/k,n-k}\,u=(-1)^{n-k}\,q^{-k}\,y_{n/k-1,n+1-k}\,T_cu
$$
where $c=(n+1)\cvar k$. Now we take $\pi=(n+k)\cvar(n+1)$ and apply $T_\pi$ to 
both sides of the last equality. As $T_\pi$ lies in the subalgebra of 
$\calH_{n+k}$ generated by the set $\{T_i\mid n<i<n+k\}$, this element commutes 
with all elements of $\calH_n$. In particular, $T_\pi$ commutes with 
$y_{n/r,n-r}$ for $r\le n$. We get 
$$
T_\pi\,y_{n/k,n-k}\,u=(-1)^{k-1}y_{n/k,n-k}\,u
$$
since $T_\pi u=(-1)^{k-1}u$, and 
$$
T_\pi\,y_{n/k-1,n+1-k}\,T_cu=y_{n/k-1,n+1-k}\,T_\pi T_cu
=y_{n/k-1,n+1-k}\,T_{(n+k)\cvar k}\,u
$$
since $\,T_\pi T_c=T_{\pi c}=T_{(n+k)\cvar k}\,$. Thus
$$
y_{n/k,n-k}\,u=(-1)^{n-1}\,q^{-k}\,y_{n/k-1,n+1-k}\,T_{(n+k)\cvar k}\,u.
$$
If $k=1$, then $\rho={(n+k)\cvar k}$, and we get the desired formula. For $k>1$ 
we use induction on $k$. Note that $\,T_{(n+k)\cvar k}\,u\in\Up^{(k-1,n)}\ot V\,$ 
by Lemma 2.1. By the induction hypothesis we may assume that
$$
y_{n/k-1,n+1-k}\,T_{(n+k)\cvar k}\,u
=(-1)^{(k-1)(n-1)}\,q^{-(k-1)k/2}\,T_{\rho'}T_{(n+k)\cvar k}\,u
$$
where $\rho'$ is the longest $k-1,n$-shuffle. Lemma 1.3 shows that 
$\rho=\rho'\cdot\bigl((n+k)\cvar k\bigr)$ and 
$T_\rho=T_{\rho'}T_{(n+k)\cvar k}$. This yields the conclusion.

If $k=n$, then $\frS_{k,n-k}=\frS_n$, and so $\,y_{n/k,n-k}=1$.
\endproof

\proclaim
Lemma 2.4.
Suppose that $\dim\Up^{(n)}=1$ and $\Up^{(n+1)}=0$ for some $n>0$. Fix 
a nonzero element $t\in\Up^{(n)}$. For each $k=0,\ldots,n$ the bilinear 
pairing
$$
\be_k:\Up^{(k)}\times\Up^{(n-k)}\to\bbk
$$
defined by the rule $\,\be_k(u,w)\,t=u\star w=y_{n/k,n-k}(uw)\,$ for 
$u\in\Up^{(k)}$ and $w\in\Up^{(n-k)}$ is nondegenerate.

There are linear bases $u_1,\ldots,u_d$ for $\Up^{(k)}$ and $w_1,\ldots,w_d$ 
for $\Up^{(n-k)}$ such that
$$
t=w_1u_1+\cdots+w_du_d\,.
$$
\endproclaim

\Proof.
Since $\Up^{(n)}\sbs\Up^{(n-k,k)}$, there is an expression 
$t=w_1u_1+\cdots+w_du_d$ for some linearly independent elements 
$u_1,\ldots,u_d\in\Up^{(k)}$ and $w_1,\ldots,w_d\in\Up^{(n-k)}$. Suppose that 
$a\in\Up^{(k)}$ is such that $\be_k(a,w_i)=0$ for all $i=1,\ldots,d$. Then
$$
y_{n/k,n-k}(at)=\sum y_{n/k,n-k}(aw_i)\,u_i=\sum\be_k(a,w_i)\,tu_i=0.
$$
On the other hand, $\,y_{n/k,n-k}(at)=(-1)^{kn-k}\,q^{-k(k+1)/2}\,T_\rho(at)\,$ 
by Lemma 2.3 since $at\in\Up^{(k,n)}$. Hence $T_\rho(at)=0$, but this entails 
$a=0$ since $T_\rho$ is an invertible element of $\calH_{k+n}$. It follows 
that $\be_k$ has zero left radical, and also
$$
\dim\Up^{(k)}\le d\le\dim\Up^{(n-k)}.
$$
Replacing $k$ with $n-k$, we deduce that $\,\dim\Up^{(n-k)}\le\dim\Up^{(k)}\,$ 
as well. Hence we must have equalities in the displayed line above. This shows 
that $\be_k$ is nondegenerate and the chosen linearly independent collections 
of elements form bases for $\Up^{(k)}$ and $\Up^{(n-k)}$.
\endproof

\proclaim
Theorem 2.5.
If\/ $\dim\Up^{(n)}=1$ and\/ $\Up^{(n+1)}=0$ then the algebras\/ $\Up(V,R)$ 
and\/ $\bbS(V,R)^!$ are Frobenius. Vanishing of\/ $\Up^{(n+1)}$ is actually 
a consequence of the condition\/ $\dim\Up^{(n)}=1$ provided that\/ $\dim V>1$.
\endproclaim

\Proof.
A connected graded algebra $A=\bbk\oplus A_1\oplus A_2\cdots$ is Frobenius if 
and only if there exists $n\ge0$ such that $\dim A_n=1$, $A_{n+1}=0$ and all 
multiplication maps
$$
A_k\otimes A_{n-k}\to A_n
$$
are nondegenerate in the sense that no nonzero element of $A_k$ annihilates 
$A_{n-k}$ on the left or on the right. For the algebra $\Up(V,R)$ this 
nondegeneracy condition is equivalent to the first conclusion of Lemma 2.4. 
Similarly, the second conclusion of Lemma 2.4 implies that $\bbS(V,R)^!$ is 
Frobenius. Indeed, the multiplication maps
$$
(\Up^{(l)})\strut^*\otimes(\Up^{(k)})\strut^*\to(\Up^{(l+k)})\strut^*
$$
in the algebra $\bbS(V,R)^!$ are dual to the inclusion maps 
$\Up^{(l+k)}\hrar\Up^{(l)}\otimes\Up^{(k)}$. Take $l=n-k$. If 
$\xi\in(\Up^{(k)})\strut^*$ is such that $(\Up^{(n-k)})\strut^*\,\xi=0$, then 
$$
\sum\eta(w_i)\xi(u_i)=(\eta\otimes\xi)(t)=0
$$
for all $\eta\in(\Up^{(n-k)})\strut^*$. We must have $\xi(u_i)=0$ for all $i$, 
whence $\xi=0$.

\smallskip
Suppose now that $\dim\Up^{(n)}=1$, but $\Up^{(n+1)}\ne0$. We will show that 
this is possible only when $\dim V=1$. Let $0\ne t\in\Up^{(n)}$. Since 
$\Up^{(n+1)}=\Up^{(n)}V\cap V\Up^{(n)}$, each element of $\Up^{(n+1)}$ can 
be written as $ta$ and as $bt$ for some vectors $a,b\in V$. But for $a\ne0$  
the equality $ta=bt$ in the tensor algebra $\bbT(V)$ can only be satisfied 
when $t$ is a scalar multiple of $a^n$. It follows that $\Up^{(n)}$ and 
$\Up^{(n+1)}$ are spanned, respectively, by the tensors $a^n$ and $a^{n+1}$ 
for some $0\ne a\in V$.

Since $a^n\in\Up^{(n)}\sbs\Up^{(1,n-1)}$, we must have $a^{n-1}\in\Up^{(n-1)}$. 
Now $va^{n-1}\in\Up^{(1,n-1)}$ and $va^n\in\Up^{(1,n)}$ for all $v\in V$. By 
Lemma 2.1
$$
y_{n/1,n-1}(va^{n-1})\in\Up^{(n)}=\bbk a^n\qquad\hbox{and}\qquad
y_{n+1/1,n}(va^n)\in\Up^{(n+1)}=\bbk a^{n+1}.
$$
Recall that $\,y_{n+1/1,n}=qy_{n/1,n-1}+(-1)^n\,T_c\,$ where $c=(n+1)\cvar1$ 
(see Lemma 1.4). Hence
$$
T_c(va^n)=(-1)^ny_{n+1/1,n}(va^n)-(-1)^nq\,y_{n/1,n-1}(va^{n-1})a\in\Up^{(n+1)}.
$$
Since $\Up^{(n+1)}$ is an $\calH_{n+1}$-submodule of $\bbT_{n+1}(V)$ and $T_c$ 
is an invertible element of $\calH_{n+1}$, we deduce that
$$
va^n\in T_c^{-1}\Up^{(n+1)}=\Up^{(n+1)}=\bbk a^{n+1}.
$$
Hence all vectors $v\in V$ are scalar multiples of $a$, i.e., $V=\bbk a$.
\endproof

\setitemsize(iii)
\proclaim
Proposition 2.6.
Suppose that\/ $\dim\Up^{(n)}=1$ and $\Up^{(n+1)}=0$. Then the following 
conditions are equivalent:

\item(i)
the action of $y_n$ on $V^{\ot n}$ is nonzero,

\item(ii)
the algebra $\Up(V,R)$ is generated by its component $\Up^{(1)}=V,$

\item(iii)
the algebra $\Up(V,R)$ is isomorphic to a factor algebra of $\La(V,R)$.

\endproclaim

\Proof.
Clearly assertion (ii) is equivalent to surjectivity of the homomorphism 
$\bbT(V)\to\Up(V,R)$ extending the identity map $V\to V$. This homomorphism 
was described in Lemma 2.2. Its surjectivity means that 
$y_kV^{\ot k}=\Up^{(k)}$ for all $k$. In particular, (ii)$\Rar$(i).

Conversely, suppose that (i) holds. Then $y_nV^{\ot n}$ is a nonzero subspace 
of the onedimensional space $\Up^{(n)}$. Hence $y_nV^{\ot n}=\Up^{(n)}$. For 
$k<n$ put $W_k=y_kV^{\ot k}$, which is a subspace of $\Up^{(k)}$. Since 
$y_n=y_{n-k,k}\,y_{n-k,k\backslash n}$, we get
$$
\Up^{(n)}=y_nV^{\ot n}\sbs y_{n-k,k}V^{\ot n}=
y_{n-k}V^{\ot(n-k)}\ot y_kV^{\ot k}=W_{n-k}\ot W_k\,,
$$
and so $t=\sum w_iu_i\in W_{n-k}\ot W_k$, in the notation of Lemma 2.4. It 
follows that $u_i\in W_k$ for all $i$. Hence $W_k=\Up^{(k)}$ since $\Up^{(k)}$ 
is spanned by $\{u_i\}$. This shows that (i)$\Rar$(ii).

Next, (ii)$\Rar$(iii) since the homomorphism of algebras $\bbT(V)\to\Up(V,R)$ 
factors through $\La(V,R)$. Finally, (iii)$\Rar$(ii) since the algebra 
$\La(V,R)$ is generated by its component of degree 1.
\endproof

It can be seen that conditions (i)--(iii) of Proposition 2.6 are satisfied at 
least when $n$ is small or under some restrictions on $q$. For a later 
reference we record here just one easy conclusion:

\proclaim
Corollary 2.7.
Suppose that\/ $\dim V>1$ and\/ $\dim\Up^{(3)}=1$. Then the action of $y_3$ 
on $V^{\ot3}$ is nonzero.
\endproclaim

\Proof.
We have seen in Theorem 2.5 that $\Up^{(4)}=0$. Note that $y_{2/1,1}=y_2=q-T_1$ and $(q-T_1)V^{\ot2}=\Im(q\cdot\Id-R)=\Up^{(2)}$. 
This means that $V\star V=\Up^{(2)}$. Also, $V\star\Up^{(2)}=\Up^{(3)}$ by Lemma 
2.4 applied with $n=3$. Thus the algebra $\Up(V,R)$ is generated by $V$.
\endproof

\proclaim
Theorem 2.8.
Suppose that\/ $\dim V>1$ and\/ $\dim\La_n(V,R)=1$ for some $n>0$. Then 
$\La_i(V,R)=0$ for all $i>n$ and the algebra\/ $\La(V,R)$ is Frobenius.
\endproclaim

\Proof.
Since $\,\La(V,R)\cong\bbS(V^*,R^*)^!$, the desired conclusion is just a 
reformulation of Theorem 2.5 with $R$ changed to $R^*$.
\endproof

If $[n]!_q\ne0$, then for all $k\le n$ the map $\La_k(V,R)\to\Up^{(k)}$ 
induced by the action of $y_k$ on $V^{\ot k}$ is bijective. In this case 
$\dim\La_n(V,R)=\dim\Up^{(n)}$, and if this dimension is equal to 1, then 
$\La(V,R)$ is a Frobenius algebra isomorphic to $\Up(V,R)$.

\proclaim
Question 2.9.
Is it true without any restriction on $q$ that $\La(V,R)\cong\Up(V,R)$ under 
the assumptions of either Theorem 2.5 or Theorem 2.8?
\endproclaim

\Remark.
Let $R$ be a braiding on a vector space $V$ not necessarily satisfying the 
quadratic Hecke relation, and let $q\in\bbk$ be any eigenvalue of $R$. With 
$T_1,T_2,\ldots$ interpreted as generators of Artin's braid groups 
$B_1,B_2,\ldots$ the definitions of algebras $\bbS(V,R)$, $\La(V,R)$, 
$\Up(V,R)$ still make sense. However, the claim that $T_\pi u=(-1)^{k-1}u$ in 
the proof of Lemma 2.3 does require the equalities $(T_i-q)(T_i+1)=0$. In the 
case of an arbitrary braiding this proof breaks down, and we are unable to 
deduce further results.

Replacing $R$ with the braiding $-q^{-1}R$, we may assume that $q=-1$, and then 
$y_n$ becomes the usual symmetrizer in the group algebra of $B_n$. Therefore 
$\La(V,R)$ is a pre-Nichols algebra, while $\Up(V,R)$ is a post-Nichols 
algebra in the sense of \cite{And17, Definition 4.8}. The corresponding 
Nichols algebra is the subalgebra of $\Up(V,R)$ generated by the degree 1 
component $V$.
\endremark

\section
3. Three fundamental operators

Let $R$ be a Hecke symmetry with parameter $q$ on a finite dimensional vector 
space $V$ over a field $\bbk$. We will be concerned with certain linear 
operators $\th,\ph,\psi\in\GL(V)$ defined in terms of $R$ in the case when 
the two algebras $\Up(V,R)$ and $\bbS(V,R)^!$ are Frobenius. By Theorem 2.5 
this happens precisely when
$$
\dim\Up^{(n)}=1\qquad\hbox{and}\qquad\Up^{(n+1)}=0
$$
for some $n>0$. We make this assumption throughout the whole section, and we 
fix a nonzero tensor $t\in\Up^{(n)}\sbs V^{\ot n}$.

In general, for a Frobenius connected graded algebra $A=\bigoplus_{k\ge0}A_k$ 
with the last nonzero component $A_n$ there is a distinguished choice of 
Nakayama automorphism. This automorphism $\nu$ preserves the grading and is 
characterized by the property
$$
ba=\nu(a)b\quad\hbox{whenever $a\in A_k$ and $b\in A_{n-k}$ for some $k$}.
$$
It commutes with every automorphism of $A$ which preserves the 
grading of $A$.

If $A$ is generated by its degree 1 component $A_1=V$, then $A\cong\bbT(V)/I$ 
where $I$ is a graded ideal of $\bbT(V)$. Since $I_n$ has codimension 1 in 
$\bbT_n(V)$, there is a linear function $f\in\bbT_n(V)^*$ such that $\,\Ker 
f=I_n$. The characteristic property of $\nu$ translates into the twisted 
cyclicity property for the function $f$:
$$
f(wv)=f\bigl(\ph(v)w),\qquad v\in V,\ w\in\bbT_{n-1}(V),
$$
where $\ph:V\to V$ is the restriction of $\nu$ to $A_1$. Nondegeneracy of the 
multiplication map $\,A_k\times A_{n-k}\to A_n\,$ implies that the 
$k\mskip1mu$th homogeneous component $I_k$ of the ideal $I$ coincides with the 
left radical of the bilinear pairing $\,\bbT_k(V)\times\bbT_{n-k}(V)\to\bbk$ 
defined by the rule $(u,w)\mapsto f(uw)$. Thus such an algebra $A$ is completely 
determined by the function $f$.

This description of graded Frobenius algebras generated in degree 1 by means 
of functions is well known (see, e.g., \cite{Bond-P93, Prop. 8.1}). In a more 
recent terminology one refers to $f$ as a \emph{preregular form} or, regarding 
$f$ as an element of $\bbT_n(V^*)$, a \emph{twisted potential} (see 
\cite{Dub14}).

Further on in this section we will denote by $\nu$ the Nakayama automorphism 
of the algebra $\Up(V,R)$ and by $\ph$ its restriction to the degree 1 component 
$\Up^{(1)}=V$. In terms of the pairings introduced in Lemma 2.4 we have
$$
\be_{n-k}(b,a)=\be_k\bigl(\nu(a),b\bigr),\qquad 
a\in\Up^{(k)},\ \,b\in\Up^{(n-k)}.
$$
It will be seen in the proof of Theorem 3.8 that the action of $\nu$ on the 
$k\mskip1mu$th homogeneous component $\Up^{(k)}$ is given by the linear operator 
$\ph^{\ot k}\in\GL(V^{\ot k})$, but this is not clear right now. In the case 
when $\Up(V,R)$ is generated by $V$, we can associate to the algebra $\Up(V,R)$ 
a linear function $f\in\bbT_n(V)^*$ with the twisted cyclicity property (see 
Proposition 3.12).

The graded Frobenius algebra $\bbS(V,R)^!$ is generated by its degree 1 
component $V^*$. Therefore its Nakayama automorphism is determined by action 
on elements of degree 1. We denote by $\psi:V\to V$ the linear operator such 
that the dual operator $\psi^*:V^*\to V^*$ is the restriction of the Nakayama 
automorphism of $\bbS(V,R)^!$.

The linear function $\bbT_n(V^*)\to\bbk$ associated with $\bbS(V,R)^!$ may be 
identified with the tensor $t$, and the twisted cyclicity property of $t$ means 
that whenever $\,t=\sum x_it_i\,$ for some elements $\,x_i\in V\,$ and 
$\,t_i\in\bbT_{n-1}(V)$ we also have
$$
t=\sum t_i\psi(x_i).
$$
It has been seen in Lemma 2.4 that $t$ can be expressed as $\,\sum x_it_i\,$ 
with $\{x_i\}$ being a linear basis for $V$ and $\{t_i\}$ a linear basis for 
$\Up^{(n-1)}$.

\proclaim
Proposition 3.1.
There are linear operators $\th,\thov\in\GL(V)$ such that
$$
T_{(n+1)\cvar1}(vt)=t\,\th(v),\qquad
T_{1\cvar(n+1)}(tv)=\thov(v)\,t
$$
for all $v\in V$. Moreover, $\,\thov=q^{n+1}\th^{-1}$.
\endproclaim

\Proof.
By Lemma 2.1 $\,T_{(n+1)\cvar1}$ maps $\Up^{(1,n)}=V\ot\Up^{(n)}$ into 
$\Up^{(n,1)}=\Up^{(n)}\ot V$. Similarly, $T_{1\cvar(n+1)}$ maps $\Up^{(n,1)}$ 
into $\Up^{(1,n)}$. Since the 1-dimensional space $\Up^{(n)}$ is spanned by 
$t$, we have $\Up^{(1,n)}=Vt$ and $\Up^{(n,1)}=t\mskip1mu V$. Hence $\th$ and 
$\thov$ are well defined. We will check that $\,\th\circ\thov=q^{n+1}\Id$.

Lemma 2.1 shows also that $\,y_{n+1/n}\Up^{(n,1)}\sbs\Up^{(n+1)}=0$. 
Since
$$
y_{n+1/n}=\sum_{i=1}^{n+1}\,(-1)^{n+1-i}\,q^{i-1}\,T_{i\cvar(n+1)},
$$
it follows that\quad
$T_{1\cvar(n+1)}u=\sum\limits_{i=2}^{n+1}\,(-1)^i\,q^{i-1}\,T_{i\cvar(n+1)}u$,\quad 
and therefore
$$
T_{(n+1)\cvar1}T_{1\cvar(n+1)}u=\sum_{i=2}^{n+1}\,(-1)^i\,q^{i-1}\,
T_{(n+1)\cvar1}T_{i\cvar(n+1)}u
$$
for all $u\in\Up^{(n,1)}$.
For $2\le i\le n+1$ easy computations in the group $\frS_{n+1}$ yield
$$
\eqalign{
\bigl((n+1)\cvar1\bigr)\bigl(i\cvar(n+1)\bigr)
&{}=\bigl((i-1)\cvar n\bigr)\bigl((n+1)\cvar1\bigr)\cr
&{}=\bigl((i-1)\cvar n\bigr)\bigl((n+1)\cvar n\bigr)\bigl(n\cvar1\bigr)\cr
&{}=\bigl((i-1)\cvar(n+1)\bigr)\bigl(n\cvar1\bigr).
}
$$
Note also that
$$
T_{(n+1)\cvar1}T_{i\cvar(n+1)}=T_{((n+1)\cvar1)(i\cvar(n+1))}
=T_{(i-1)\cvar(n+1)}T_{n\cvar1}
$$
and $\,T_{n\cvar1}u=(-1)^{n-1}u\,$ for all $u\in\Up^{(n,1)}$. Hence
$$
T_{(n+1)\cvar1}T_{1\cvar(n+1)}u=\sum_{i=1}^n\,(-1)^{n-i}q^i\,T_{i\cvar(n+1)}u
=q^{n+1}u-q\,y_{n+1/n}\,u=q^{n+1}u.
$$
Taking $u=tv$, we get $\th\bigl(\thov(v)\bigr)=q^{n+1}v$, as claimed.
\endproof

\proclaim
Lemma 3.2.
Let $\rho\in\frS_{k+n}$ be the longest $k,n$-shuffle where $k$ is any positive 
integer. Then
$$
T_\rho(ut)=t\,\th^{\ot k}(u),\qquad T_{\rho^{-1}}(tu)=\thov^{\ot k}(u)\,t
$$
for all $\,u\in V^{\ot k}$.
\endproclaim

\Proof.
If $k=1$, then $\rho=(n+1)\cvar1$, and the formulas in the statement of Lemma 
3.2 are those used in the definition of $\th$ and $\thov$. Let $k>1$. It suffices 
to prove the formulas for $u=av$ where $a\in V^{\ot(k-1)}$ and $v\in V$. Let 
$\rho'\in\frS_{k+n-1}$ be the longest $k-1,n$-shuffle. Proceeding by induction 
on $k$ we may assume that
$$
T_{\rho'}(at)=t\,\th^{\ot(k-1)}(a),\qquad 
T_{\rho'^{-1}}(ta)=\thov^{\ot(k-1)}(a)\,t.
$$
By Lemma 1.3 $\rho=\rho'c$ and $T_\rho=T_{\rho'}T_c$ where $c=(n+k)\cvar k$. 
Then $\rho^{-1}=c^{-1}\rho'^{-1}$ where $c^{-1}=k\cvar(n+k)$. Since 
$\ell(\rho)=\ell(\rho')+\ell(c)$ and the lengths of permutations do not change 
after taking inverses, we also have $T_{\rho^{-1}}=T_{c^{-1}}T_{\rho'^{-1}}$. 
Noting that $T_c$ and $T_{c^{-1}}$ act on $V^{\ot(k+n)}$, respectively, as 
$$
\Id_V^{\ot(k-1)}\ot\ T_{(n+1)\cvar 1}\quad\hbox{and}\quad
\Id_V^{\ot(k-1)}\ot\ T_{1\cvar(n+1)},
$$
we get\quad $\displaystyle T_\rho(ut)=T_{\rho'}T_c(avt)
=T_{\rho'}\bigl(at\th(v)\bigr)=t\,\th^{\ot(k-1)}(a)\th(v)=t\,\th^{\ot k}(u)$\quad 
and
$$
T_{\rho^{-1}}(tu)=T_{c^{-1}}T_{\rho'^{-1}}(tav)
=T_{c^{-1}}\bigl(\thov^{\ot(k-1)}(a)tv\bigr)
=\thov^{\ot(k-1)}(a)\thov(v)t=\thov^{\ot k}(u)\,t,
$$
as claimed.
\endproof

\proclaim
Corollary 3.3.
The linear operator $\th\ot\th$ on $V^{\ot2}$ commutes with $R$. In particular, 
$\th$ extends to automorphisms of $\La(V,R)$ and $\bbS(V,R)$.
\endproclaim

\Proof.
Lemma 3.2 with $k=2$ shows that $\,T_\rho(ut)=t\,(\th\ot\th)(u)\,$ for all 
$\,u\in V^{\ot2}$. By Lemma 1.3 $\,T_\rho T_1=T_{n+1}T_\rho\,$. Hence
$$
T_\rho(R(u)t)=T_\rho T_1(ut)=T_{n+1}T_\rho(ut)
=T_{n+1}\bigl(t\mskip1mu(\th\ot\th)(u)\bigr)=t\,\bigl(R(\th\ot\th)(u)\bigr),
$$
which implies that $(\th\ot\th)R(u)=R(\th\ot\th)(u)$ for all $u\in V^{\ot2}$. 
As a consequence, the ideals of the tensor algebra $\bbT(V)$ defining its 
factor algebras $\bbS(V,R)$ and $\La(V,R)$ are stable under the automorphism 
of $\bbT(V)$ extending $\th$.
\endproof

\proclaim
Corollary 3.4.
We have $\,\th^{\ot n}(t)=q^{n(n+1)/2}\,t$.
\endproclaim

\Proof.
Take $k=n$ and $u=t\in\Up^{(n)}$ in Lemma 3.2. Since $ut=t^2\in\Up^{(n,n)}$, 
we know from Lemma 2.3 that $T_\rho(ut)=q^{\mskip2mu n(n+1)/2}\,t^2$. This 
yields the formula for $\th^{\ot n}(t)$.
\endproof

We want to use both $\th$ and $\thov$. For this we will need also a 
right hand version of Lemma 2.3:

\proclaim
Lemma 3.5.
Let $1\le k\le n,$ and let $\rho\in\frS_{k+n}$ be the longest $k,n$-shuffle. Then
$$
y_{n/n-k,k}^{(k)}u=(-1)^{kn-k}\,q^{-k(k+1)/2}\,T_{\rho^{-1}}u
\quad\hbox{for all $\,u\in\Up^{(n,k)}$}
$$
where $y_{n/n-k,k}^{(k)}\in\calH_{k+n}$ is the $k$ step shift of 
$y_{n/n-k,k}$.
\endproclaim

\Proof.
We follow the proof of Lemma 2.3. Here we use the identity
$$
y_{n+1/n+1-k,k}=q^ky_{n/n-k,k}^{(1)}+(-1)^{n+1-k}y_{n/n+1-k,k-1}^{(1)}T_{1\cvar(n+2-k)}
$$
in $\calH_{n+1}$ (cf. Lemma 1.4). After taking the $k-1$ step shifts of 
elements it gives
$$
y_{n/n-k,k}^{(k)}u=(-1)^{n-k}\,q^{-k}y_{n/n+1-k,k-1}^{(k)}T_{k\cvar(n+1)}u
$$
since $y_{n+1/n+1-k,k}^{(k-1)}u=0$. Applying $T_{1\cvar k}$ to both sides, 
we get
$$
y_{n/n-k,k}^{(k)}u=(-1)^{n-1}\,q^{-k}y_{n/n+1-k,k-1}^{(k)}T_{1\cvar(n+1)}u.
$$
For $k=1$ this is the desired equality with $\rho^{-1}=1\cvar(n+1)$. For $k>1$ 
we use induction on $k$ noting that $T_{1\cvar(n+1)}u\in V\ot\Up^{(n,k-1)}$
and
$$
T_{\rho^{-1}}=T_{k\cvar(n+k)}\cdots T_{2\cvar(n+2)}T_{1\cvar(n+1)}=
T_{(\rho')^{-1}}^{\,(1)}\,T_{1\cvar(n+1)}
$$
where $\rho'\in\frS_{k+n-1}$ is the longest $k-1,n$-shuffle.
\endproof

\proclaim
Proposition 3.6.
Given $1\le k\le n,$ let $u_i\in\Up^{(k)},$ $\,w_i\in\Up^{(n-k)}$ 
be elements such that
$$
t=\sum u_iw_i=\sum w_i\,\psi^{\ot k}(u_i).
$$
Then
$$
\eqalign{
\th^{\ot k}(a)&{}=(-1)^{kn-k}\,q^{k(k+1)/2}\,\sum\be_k(a,w_i)\,\psi^{\ot k}(u_i),\cr
\thov^{\ot k}(a)&{}=(-1)^{kn-k}\,q^{k(k+1)/2}\,\sum\be_{n-k}(w_i,a)\,u_i
}
$$
for all $a\in\Up^{(k)}$. In particular,
$$
\th(v)=(-1)^{n-1}q\sum\be_1(v,t_i)\psi(x_i),\qquad
\thov(v)=(-1)^{n-1}q\sum\be_{n-1}(t_i,v)x_i
$$
for all $v\in V$ where $x_i\in V$ and $t_i\in\Up^{(n-1)}$ are such that 
$\,t=\sum x_it_i=\sum t_i\psi(x_i)$.
\endproclaim

\Proof.
The second equality in the requested expression for $t$ follows from the twisted 
cyclicity of $t$, and such an expression exists since $\Up^{(n)}\sbs\Up^{(k,n-k)}$. 
By Lemmas 3.2 and 2.3
$$
t\,\th^{\ot k}(a)=T_\rho(at)=(-1)^{kn-k}\,q^{k(k+1)/2}\,y_{n/k,n-k}(at).
$$
Substituting $\,t=\sum w_i\,\psi^{\ot k}(u_i)\,$, we get
$$
y_{n/k,n-k}(at)=\sum y_{n/k,n-k}(aw_i)\,\psi^{\ot k}(u_i)
=t\cdot\sum\be_k(a,w_i)\,\psi^{\ot k}(u_i)
$$
by the definition of $\be_k$ in Lemma 2.4. The formula for $\th^{\ot k}$ 
follows after cancelling the common factor $t$. The formula for 
$\thov^{\ot k}$ is derived similarly, using Lemma 3.5 to expand 
$\,T_{\rho^{-1}}(ta)$.
\endproof

\proclaim
Corollary 3.7.
We have $\,\th^{\ot k}\bigl(\nu(a)\bigr)=(\psi\thov)^{\ot k}(a)\,$ for all 
$\,a\in\Up^{(k)}$.
\endproclaim

\Proof.
Since $\,\be_k\bigl(\nu(a),w_i\bigr)=\be_{n-k}(w_i,a)$, Proposition 3.6 yields
$$
\th^{\ot k}\bigl(\nu(a)\bigr)
=\sum(-1)^{kn-k}\,q^{k(k+1)/2}\,\be_{n-k}(w_i,a)\,\psi^{\ot k}(u_i)
=\psi^{\ot k}\bigl(\thov^{\ot k}(a)\bigr),
$$
as required.
\endproof

\proclaim
Theorem 3.8.
The operators $\ph,\psi,\th$ pairwise commute and $\,\psi=q^{-n-1}\ph\,\th^2$. 
Moreover, they lie in the subgroup of $\GL(V)$ consisting of all invertible 
linear operators $\chi:V\to V$ such that $\chi\ot\chi$ commutes with $R$.
\endproclaim

\Proof.
Recall that $\ph$ is the degree 1 component of $\nu$. Therefore the case $k=1$ 
of Corollary 3.7 gives $\,\th\ph=\psi\thov$. Since $\,\thov=q^{n+1}\th^{-1}$, 
we get $\,\psi=q^{-n-1}\th\,\ph\,\th$.

Now $\,\th^{\ot k}\ph^{\ot k}=(\th\ph)^{\ot k}=(\psi\thov)^{\ot k}\,$ for 
each $k>0$. Comparing this with the identity of Corollary 3.7, we deduce that
$$
\nu(a)=\ph^{\ot k}(a)\quad\hbox{for all $a\in\Up^{(k)}$}.
$$
Let $u,v\in V$. Since $\nu$ is an automorphism of the algebra $\Up(V,R)$, we get
$$
\ph^{\ot2}(u\star v)=\ph(u)\star\ph(v).
$$
But $u\star v=y_2(uv)=(q-T_1)(uv)$. Similarly, 
$\ph(u)\star\ph(v)=(q-T_1)\bigl(\ph(u)\ph(v)\bigr)$. Therefore the displayed 
identity above means precisely that $\ph^{\ot2}$ commutes with the operator 
$q\cdot\Id-R$. Then $\ph^{\ot2}$ commutes with $R$. In Corollary 3.4 we have 
seen that $\th^{\ot2}$ commutes with $R$. Then so too does $\psi^{\ot2}$.

The fact that $\th^{\ot2}$ commutes with $R$ implies that for each $k>0$ the 
linear operator $\th^{\ot k}$ is an $\calH_k$-module endomorphism of $V^{\ot k}$. 
As a consequence, $\Up^{(k)}$ is stable under $\th^{\ot k}$, and for 
$a\in\Up^{(k)}$ and $b\in\Up^{(l)}$ we have
$$
\th^{\ot(k+l)}(a\star b)=\th^{\ot(k+l)}\bigl(y_{k+l/k,l}(ab)\bigr)=
y_{k+l/k,l}\bigl(\th^{\ot k}(a)\,\th^{\ot l}(b)\bigr)=
\th^{\ot k}(a)\star\th^{\ot l}(b).
$$
Thus the linear operators $\,\th^{\ot k}|_{\Up^{(k)}}$, $\,k\ge0$, define an 
automorphism of the algebra $\Up(V,R)$ preserving the grading. Any such an 
automorphism commutes with the Nakayama automorphism $\nu$. In particular, 
$\,\th\ph=\ph\th$. The rest is clear.
\endproof

\proclaim
Proposition 3.9.
For $1\le k\le n$ denote by $\xi_k$ and $\eta_k,$ respectively, the 
restrictions of the linear operators $(\psi^{-1}\th)^{\ot k}$ and 
$(\psi\mskip1mu\thov)^{\ot k}$ to the subspace $\Up^{(k)}\sbs V^{\ot k}$. Then
$$
\tr\xi_k=\tr\eta_k=(-1)^{kn-k}\,q^{k(k+1)/2}\left[{n\atop k}\right]_q
$$
where $\left[{n\atop k}\right]_q$ is the $q$-binomial coefficient. 
In particular,
$$
\tr\psi^{-1}\th=\tr\psi\mskip1mu\thov=(-1)^{n-1}q\,[n]_q\,.
$$
\endproclaim

\Proof.
Put $C=(-1)^{kn-k}\,q^{k(k+1)/2}$. In the notation of Proposition 3.6 we have
$$
\eqalign{
\xi_k(a)&{}=(\psi^{-1})^{\ot k}\,\th^{\ot k}(a)
=C\,\sum\be_k(a,w_i)\,u_i,\cr 
\eta_k(a)&{}=\psi^{\ot k}\,\thov^{\ot k}(a)
=C\,\sum\be_{n-k}(w_i,a)\,\psi^{\ot k}(u_i)
}
$$
for $a\in\Up^{(k)}$. Note that for each $u\in\Up^{(k)}$ and 
$g\in(\Up^{(k)})^*$ the linear operator on $\Up^{(k)}$ defined by the 
rule $a\mapsto g(a)u$, $\,a\in\Up^{(k)}$, has trace $g(u)$. This yields
$$
\tr\xi_k=C\,\sum\be_k(u_i,w_i),\qquad
\tr\eta_k=C\,\sum\be_{n-k}\bigl(w_i,\psi^{\ot k}(u_i)\bigr),
$$
i.e., by the definition of bilinear pairings in Lemma 2.4,
$$
\eqalign{
(\tr\xi_k)\,t&{}=C\,y_{n/k,n-k}\bigl({\textstyle\sum}u_iw_i\bigr)
=C\,y_{n/k,n-k}\mskip1mu t,\cr
(\tr\eta_k)\,t&{}=C\,y_{n/n-k,k}\bigl({\textstyle\sum}w_i\psi^{\ot k}(u_i)\bigr)
=C\,y_{n/n-k,k}\mskip1mu t.
}
$$
Since $T_\si\mskip1mu t=(-1)^{\ell(\si)}$ for each $\si\in\frS_n$, we have
$$
y_{n/k,n-k}\mskip1mu t
=\sum_{\si\in\calD(\frS_n/\frS_{k,n-k})}q^{k(n-k)-\ell(\si)}\,t
=\sum_{\si\in\calD(\frS_n/\frS_{k,n-k})}q^{\ell(\si)}t
=\left[{n\atop k}\right]_qt.
$$
Evaluation of $y_{n/n-k,k}\mskip1mu t$ gives the same result. Hence 
$\,\tr\xi_k=\tr\eta_k=C\left[{n\atop k}\right]_q$.
\endproof

\proclaim
Corollary 3.10.
Let $\xi_k$ be defined as in Proposition 3.9. Then
$$
\tr\xi_{n-k}=q^{(n-2k)(n+1)/2}\tr\xi_k.
$$
\endproclaim

\Proof.
This follows from the formulas for traces in Proposition 3.9 in view of the 
equality $\,\left[{n\atop n-k}\right]_q=\left[{n\atop k}\right]_q$.
\endproof

\proclaim
Corollary 3.11.
Suppose that $n$ is odd, say $n=2m+1,$ and $\th=\la\psi$ where $\la\in\bbk$. 
If the characteristic of\/ $\bbk$ does not divide the dimension 
of\/ $\Up^{(m)},$ then $\la=q^{m+1}$.
\endproclaim

\Proof.
Since $\psi^{-1}\th=\la\cdot\Id$, we have $\,\tr\xi_k=\la^k\dim\Up^{(k)}\,$ 
for each $k$. Hence for $k=m$ the equality of Corollary 3.10 is written as
$$
\la^{m+1}\dim\Up^{(m+1)}=q^{m+1}\la^m\dim\Up^{(m)}.
$$
Since the bilinear pairing $\be_m:\Up^{(m)}\times\Up^{(m+1)}\to\bbk$ is 
nondegenerate, the two spaces here have equal dimension, and the conclusion follows.
\endproof

\proclaim
Proposition 3.12.
Suppose that $y_nV^{\ot n}\ne0$. Let $f\in\bbT_n(V)^*$ be a linear function 
such that
$$
\Ker f=\{u\in\bbT_n(V)\mid y_nu=0\}.
$$
Put $\,I'=\bigoplus I'_k\,$ where
$$
I'_k=\{u\in\bbT_k(V)\mid f(uw)=0\hbox{ for all }w\in\bbT_{n-k}(V)\},\qquad 
k=0,\ldots,n.
$$
Then $f$ satisfies the twisted cyclicity condition
$$
f(wv)=f\bigl(\ph(v)w),\qquad v\in V,\ w\in\bbT_{n-1}(V),
$$
and $I'$ is a graded ideal of\/ $\bbT(V)$ such that $I'_0=I'_1=0,$ $\,I'_2$ 
coincides with the component $I_2$ of the ideal $I$ defining the algebra\/ 
$\La(V,R),$ and\/ $\bbT(V)/I'$ is a Frobenius algebra isomorphic to\/ 
$\Up(V,R)$. If\/ $[n]!_q\ne0,$ then $I'=I$.
\endproclaim

\Proof.
By Proposition 2.6 the algebra $\Up(V,R)$ is generated by $\Up^{(1)}=V$. Hence
$$
\Up(V,R)\cong\bbT(V)/K
$$
where $K$ is the kernel of the homomorphism $\,\bbT(V)\to\Up(V,R)\,$ described 
in Lemma 2.2. This homomorphism is defined by the actions of the antisymmetrizers 
$y_k$. In particular, $\,\Ker f$ coincides with the $n\mskip1mu$th homogeneous 
component of the ideal $K$. So it follows that $f$ is exactly the linear function, 
unique up to a scalar multiple, associated with the algebra $\Up(V,R)$, as 
discussed at the beginning of this section. It satisfies the twisted 
cyclicity, as stated. Since $I'_k$ is the left radical of the pairing
$$ 
\bbT_k(V)\times\bbT_{n-k}(V)\to\bbk,\qquad(u,w)\mapsto f(uw),
$$
we deduce that $I'=K$. Since $y_0=y_1=1$, we get $I'_0=I'_1=0$. By Lemma 2.2 
$\,I\sbs I'$. We have $I'_2=\{u\in V^{\ot2}\mid y_2u=0\}=I_2$ since $y_2$ acts 
on $V^{\ot2}$ by means of the linear operator $q\cdot\Id-R$. 

Suppose that $[n]!_q\ne0$. Then $I'_k=I_k$ for all $k\le n$ by Lemma 1.1. In 
particular, $\,\dim\La_n(V,R)=1$, and then $I_k=\bbT_k(V)$ for all $k>n$, 
as is seen in Theorem 2.8 (in the trivial case $\dim V=1$ the equality 
$\Up^{(n+1)}=0$ is possible only when $n=1$ and $I_2=I'_2=\bbT_2(V)$). Hence 
$I'_k=I_k$ for all $k$.
\endproof

Under a restriction on $q$ the formulas of Proposition 3.6 can 
be rewritten in terms of the function $f$. We thus obtain a coordinate-free 
interpretation of Gurevich's commutation formulas \cite{Gur90, Prop. 5.7}. In 
the next proposition we use a normalization of $f$ which permits us to 
consider the case $[n]_q=0$. The coefficients in \cite{Gur90} correspond to a 
different normalization condition $f(t)=1$ which requires $[n]!_q\ne0$.

\proclaim
Proposition 3.13.
Suppose that $[n-1]!_q\ne0$. Let $f\in\bbT_n(V)^*$ be the linear function 
defined by the rule
$$
y_nu=[n-1]!_q\,f(u)t\quad\hbox{ for $u\in\bbT_n(V)$}.
$$
Then $\,f(t)=[n]_q\,$ and $f$ satisfies all conclusions of Proposition 3.12. 

\smallskip
Writing $\,t=\sum x_it_i\,$ for some elements $\,x_i\in V\,$ and 
$\,t_i\in\Up^{(n-1)},$ we have
$$
\th(v)=(-1)^{n-1}q\sum f(vt_i)\psi(x_i)\quad and\quad
\thov(v)=(-1)^{n-1}q\sum f(t_iv)x_i\,.
$$
for all $v\in V$.
\endproclaim

\Proof.
Let $v\in V$ and $w\in\Up^{(n-1)}$. Then $y_{n-1}w=[n-1]!_q\,w$, and therefore
$$
[n-1]!_q\,f(vw)t=y_n(vw)=v\star(y_{n-1}w)=[n-1]!_q\,v\star w
=[n-1]!_q\,\be_1(v,w)t.
$$
It follows that $f(vw)=\be_1(v,w)$. Similarly, $f(wv)=\be_{n-1}(w,v)$. Hence 
the formulas for $\th$ and $\thov$ repeat those of Proposition 3.6.

Since $\be_1\ne0$, we see that $f\ne0$. This means that $y_nV^{\ot n}\ne0$, 
and Proposition 3.12 does apply. That $\,f(t)=[n]_q\,$ is clear since 
$\,y_nt=[n]!_q\,t$.
\endproof

\Remark.
Let $\tau:V^{\ot2}\to V^{\ot2}$ be the usual flip of tensorands 
$v_1\ot v_2\mapsto v_2\ot v_1$. Then $R\op=\tau R\,\tau$ is a Hecke symmetry 
with the same parameter $q$ such that the identity map on $V$ extends to 
antiisomorphisms of algebras
$$
\La(V,R)\to\La(V,R\op),\qquad\Up(V,R)\to\Up(V,R\op)
,\qquad\bbS(V,R)\to\bbS(V,R\op).
$$
When $R$ is replaced by $R\op$ the operators $\ph,\psi,\th$ change to 
$\ph^{-1},\psi^{-1},\,\thov$.
\endremark

\section
4. Regular algebras of type A and the Hessian group

The Artin-Schelter regular graded algebras with 3 generators of degree 1 and 3 
defining relations of degree 2 were classified by Artin, Tate and Van den 
Bergh \cite{Ar-TV90} in terms of pairs $(E,\si)$ where $E$ is either the 
projective plane $\bbP^2$ or a cubic divisor in $\bbP^2$, and $\si$ is an 
automorphism of $E$. The \emph{point scheme} $E$ parametrizes so-called point 
modules for the corresponding algebra, while $\si$ can be interpreted as 
transformation of point modules by shift of degrees. In this classification 
the algebras of \emph{type $A$} are distinguished by the properties that $E$ 
is a smooth cubic curve and $\si$ is a translation by a point of the abelian 
variety $\Pic_0E\cong E$ \cite{Ar-TV90, 4.13}.

Automorphisms of the graded algebra $A$ corresponding to $(E,\si)$ are 
extensions of the linear operators $\th$ acting on the 3-dimensional component 
$V=A_1$ such that the induced transformation $\th'$ of the projective plane 
$\bbP^2=\bbP(V^*)$ leaves $E$ stable and the restriction of $\th'$ to $E$ 
commutes with $\si$ \cite{Ar-TV91, Prop. 8.8}. However, this description of 
automorphisms is not quite explicit, and moreover in the present paper we do 
not use the shift $\si$ and the whole geometric approach of \cite{Ar-TV90}.

For an algebra of type $A$ Proposition 4.1 describes automorphisms in terms of 
the Hessian group $G\sbs\PGL(V^*)$ of order 216. We refer to the book of 
Brieskorn and Kn\"orrer \cite{Br-Kn} for details on this group. In the proof 
of Theorem 5.1 it will be important that all operators $\th$ whose images 
$\th'$ lie in some conjugacy class of $G$ can be treated in the same way. This 
reduces the number of cases to consider.

In this section the base field $\bbk$ will be assumed to be algebraically 
closed of characteristic $\ne2,3$. We will be concerned with a family of graded 
quadratic algebras introduced by Artin and Schelter \cite{Ar-Sch87}. An algebra 
$A$ in this family has 3 generators $x_1,x_2,x_3$ and 3 defining relations
$$
ax_{i+1}x_{i-1}+bx_{i-1}x_{i+1}+cx_i^2=0,\qquad i=1,2,3,
$$
where $(a,b,c)\in\bbk^3\setm\{(0,0,0)\}$ is a triple of parameters and the 
indices $1,2,3$ are viewed as elements of $\bbZ/3\bbZ$. Thus $A=\bbT(V)/J$ 
where $V$ is the 3-dimensional vector space spanned by $x_1,x_2,x_3$ and $J$ 
is the ideal of $\bbT(V)$ generated by $t_1,t_2,t_3$ where
$$
t_i=ax_{i+1}x_{i-1}+bx_{i-1}x_{i+1}+cx_i^2\in\bbT_2(V).
$$
This algebra is denoted $\Skl_3(a,b,c)$ by association with Sklyanin's 
elliptic algebras of global dimension 4.

The algebra $A$ is Artin-Schelter regular precisely when its quadratic dual 
$A^!$ is a Frobenius algebra. This happens when at least 2 of the parameters 
$a,b,c$ are nonzero and the two equalities $a^3=b^3=c^3$ do not hold 
simultaneously. We will always assume that these conditions are satisfied. In 
this case $J_2V\cap VJ_2\sbs\bbT_3(V)$ is the 1-dimensional space $\<t\>$ 
spanned by the tensor
$$
t=\sum\bigl(ax_{i-1}x_ix_{i+1}+bx_{i+1}x_ix_{i-1}+cx_i^3)
=\sum x_it_i=\sum t_ix_i.
$$
The image of $t$ in the symmetric algebra $\bbS(V)\cong\bbk[x_1,x_2,x_3]$ is 
the polynomial
$$
t^S=3(a+b)x_1x_2x_3+c(x_1^3+x_2^3+x_3^3).
$$
The equation $t^S=0$ defines a curve in the projective plane $\bbP(V^*)$ 
associated with the dual space $V^*$. This curve is smooth when $c\ne0$ and 
$(a+b)^3+c^3\ne0$. The point scheme of $A$ is the curve defined by a different 
equation
$$
abc(x_1^3+x_2^3+x_3^3)=(a^3+b^3+c^3)x_1x_2x_3.
$$
The 3-dimensional Artin-Schelter regular quadratic algebras of type $A$ are 
precisely the algebras $\Skl_3(a,b,c)$ with a smooth point scheme. Here 
smoothness amounts to
$$
abc\ne0\qquad\hbox{and}\qquad(a^3+b^3+c^3)^3\ne27a^3b^3c^3.
$$
Note that
$$
\openup1\jot
\eqalign{
(a^3+b^3+c^3)^3-27a^3b^3c^3
&=\prod_{i=1}^3\prod_{j=1}^3\,\bigl(\ep^ia+\ep^jb+c\bigr)\cr
&=\bigl((a+b)^3+c^3\bigr)\bigl((a+\ep b)^3+c^3\bigr)
\bigl((a+\ep^2b)^3+c^3\bigr)
}
$$
where $\ep\in\bbk$ is a primitive cube root of 1. Hence $(a+b)^3+c^3\ne0$ 
for each algebra of type $A$, and the curve defined by the equation $t^S=0$ is 
also smooth.

Consider the linear pencil $\calL$ of cubic curves in $\bbP(V^*)$ defined 
by the equations
$$
\al x_1x_2x_3+\be(x_1^3+x_2^3+x_3^3)=0,\qquad\al,\be\in\bbk.
$$
Singular curves in $\calL$ are unions of 3 lines, while smooth ones are elliptic. 
The curves in $\calL$ have 9 common points with homogeneous coordinates
$$
\eqalign{ 
&(0:1:-1),\quad (0:1:-\ep),\quad (0:1:-\ep^2),\cr
&(-1:0:1),\quad (-\ep:0:1),\quad (-\ep^2:0:1),\cr
&(1:-1:0),\quad (1:-\ep:0),\quad (1:-\ep^2:0)
}
$$
where $\ep$ is a primitive cube root of 1, and these 9 points are precisely 
the inflection points of all nonsingular curves in $\calL$.

The Hessian group $G$ consists of all projective linear transformations of 
$\bbP(V^*)$ permuting the 9 points listed above. If $g\in G$ and $C\in\calL$, 
then $g(C)$ is a cubic passing through the same 9 points, and such a curve must 
be again a member of $\calL$. Therefore $G$ can be described as the group of 
projective linear transformations permuting the curves in $\calL$.

Suppose that $C\in\calL$ is nonsingular. Then the Hessian group $G$ contains 
each transformation $g\in\PGL(V^*)$ such that $g(C)=C$ since the set of 
inflection points of $C$ is invariant under automorphisms of the curve. 
In other words, the group
$$
G_C=\{g\in\PGL(V^*)\mid g(C)=C\}
$$
is a subgroup of $G$. The index of $G_C$ in $G$ is equal to the cardinality of 
the orbit of $C$ with respect to the action of $G$ on $\calL$. It is known that 
this number is 12 except when the $j$-invariant of $C$ is either 0 or 1728, in 
which cases there are, respectively, 4 and 6 curves in the orbit. Thus
$$
(G:G_C)=\cases{4 & if $j(C)=0$,\cr 6 & if $j(C)=1728$,\cr 12 & otherwise.}
$$

For convenience we make use of the canonical isomorphism $\GL(V^*)\cong\GL(V)$ 
which allows us to represent projective linear transformations of $\bbP(V^*)$ 
by linear operators acting on $V$.

Let $\bbF_3$ be the 3 element field. It is known that $G$ is isomorphic to the 
semidirect product of the additive group $\bbF_3^2$ by $\SL_2(\bbF_3)$. Thus 
$G$ contains a normal subgroup $T$ of order 9. It is generated by 2 
transformations represented, respectively, by the linear operators
$$
x_1\mapsto x_2,\quad x_2\mapsto x_3,\quad x_3\mapsto x_1,\quad\hbox{and}\quad
x_1\mapsto\ep x_1,\quad x_2\mapsto\ep^2x_2,\quad x_3\mapsto x_3.
$$
Each curve in the pencil $\calL$ is stable under the action of $T$, and one can 
see that smooth curves have no points fixed by a nontrivial element of $T$. 
This means that $T$ acts on each elliptic curve in $\calL$ by translations. 
Accordingly, we call $T$ the \emph{subgroup of translations} in $G$.

Now $G/T\cong\SL_2(\bbF_3)$. The group $\SL_2(\bbF_3)$ has a center of order 
2. The nontrivial central element of $G/T$ has a representative in $\GL(V)$ 
given by the assignments
$$
x_1\mapsto x_2,\quad x_2\mapsto x_1,\quad x_3\mapsto x_3.
$$
Each curve in $\calL$ is stable under this transformation as well. Denote by $Z$ 
the subgroup of $G$ containing $T$ such that $Z/T$ is the center of $G/T$. 
Then $Z$ has order 18, and $G/Z$ is a group of order 12 isomorphic to the 
alternating group on 4 letters. For each smooth curve $C\in\calL$ we get $G_C=Z$ 
except for two cases when the stabilizer of $C$ is larger:
$$
(G_C:Z)=\cases{3 & if $j(C)=0$,\cr 2 & if $j(C)=1728$.}
$$

Denote by $\Aut A$ the group of those automorphisms of $A$ which preserve the 
grading of this algebra. Since $A$ is generated by $A_1=V$, we may identify 
this group with a subgroup of $\GL(V)$.

\proclaim
Proposition 4.1.
Suppose that the base field\/ $\bbk$ is algebraically closed of characteristic 
$\ne2,3$. Let $a,b,c\in\bbk$ be such that the algebra $A=\Skl_3(a,b,c)$ is 
Artin-Schelter regular and the corresponding cubic polynomial $t^S\in\bbS(V)$ 
is irreducible. So the curve $C\sbs\bbP(V^*)$ defined by the equation $t^S=0$ 
is elliptic.

Then the image $G_A$ of the group $\Aut A$ in $\PGL(V^*)$ is a subgroup of 
the Hessian group $G$. If $a\ne b,$ then $G_A=T,$ the subgroup of translations 
in $G$. If $a=b,$ then $G_A=G_C,$ the stabilizer of $C$ in $G$. In the latter 
case $j(C)\ne0,$ and
$$
(G_A:T)=\cases{4 & if $j(C)=1728$,\cr 2 & otherwise.}
$$
\endproclaim

\Proof.
The tensor $t$ is invariant under cyclic permutations of tensorands. It can be 
written as the sum $t=t^++t^-$ of a symmetric and an alternating tensors. 
Explicitly,
$$
t^+={a+b\over2}\sum_{\pi\in\frS_3}x_{\pi 1}x_{\pi 2}x_{\pi 3}+c\sum x_i^3,
\qquad
t^-={a-b\over2}\sum_{\pi\in\frS_3}\sgn(\pi)x_{\pi 1}x_{\pi 2}x_{\pi 3}.
$$
Since the space of defining relations $J_2=\<t_1,t_2,t_3\>$ is recovered from 
the tensor $t$, a linear operator $\th\in\GL(V)$ extends to an automorphism of 
$A$ if and only if the 1-dimensional subspace $\<t\>\sbs V^{\ot3}$ spanned by 
$t$ is stable under the operator $\th^{\ot3}$, if and only if $t^+$ and $t^-$ 
lie in one eigenspace of $\th^{\ot3}$. Since $\chr\bbk\ne2,3$, the space of 
symmetric tensors in $V^{\ot3}$ is isomorphic as a $\GL(V)$-module to the 
homogeneous component $\bbS_3(V)$ of $\bbS(V)$. Hence $\th^{\ot3}(t^+)\in\<t^+\>$ 
if and only if $\th t^S\in\<t^S\>$. This condition on $\th$ means precisely 
that the curve $C$ is invariant under the respective projective transformation 
$\th'$. In other words, $\th'\in G_C$. This shows that $G_A\sbs G_C$.

Note that $t^-=0$ if and only if $a=b$. If $a=b$, then $t=t^+$. In this case 
$\<t\>$ is stable under $\th^{\ot3}$ for any linear operator $\th\in\GL(V)$ 
such that $\th'\in G_C$. Therefore $G_A=G_C$. In particular, $T\sbs G_A$. 
Since $A$ is assumed to be Artin-Schelter regular, we have $a\ne0$ and 
$a^3\ne c^3$. Recall from \cite{Ar-TV90, (10.15)} that the cubic curve defined 
by the equation $x_1^3+x_2^3+x_3^3+6\ka\,x_1x_2x_3=0$ has $j$-invariant
$$
j=-2^{12}\cdot 3^3(\ka^3-1)^3\ka^3/(8\ka^3+1)^3.
$$
So $j=0$ if and only if either $\ka=0$ or $\ka^3=1$. For the curve $C$ here we 
have $\ka=a/c$, and it follows that $j(C)\ne0$. The index $(G_C:T)$ is twice
the index $(G_C:Z)$ since $(Z:T)=2$.

Suppose now that $a\ne b$. In this case $t^-\ne0$ and $\th t^-=(\det\th)t^-$ 
for each linear operator $\th\in\GL(V)$. In order that $\<t\>$ be stable under 
$\th^{\ot3}$, it is necessary and sufficient that
$$
\th t^S=(\det\th)t^S.
$$
Thus $G_A$ consists of all elements $\th'\in G$ whose representatives in 
$\GL(V)$ satisfy this condition. It does hold for the two linear operators 
representing generators of $T$. Hence $T\sbs G_A$.

If $G_A\ne T$, then we can find $\th'\in G_A$ with the property that the image 
of $\th'$ in $G/T\cong\SL_2(\bbF_3)$ has prime order, say $p$. Since $|G/T|=24$, 
there are two possibilities. If $p=2$, then the coset $\th'T$ is the generator 
of the center $Z/T$ of $G/T$ (see the next lemma), and we must have $Z\sbs G_A$. 
However, the linear operator
$$
x_1\mapsto x_2,\quad x_2\mapsto x_1,\quad x_3\mapsto x_3
$$
representing one of such elements has $\det\th=-1$, while it fixes all 
polynomials in the 2-dimensional subspace $L\sbs\bbS_3(V)$ spanned by 
$x_1x_2x_3$ and $x_1^3+x_2^3+x_3^3$. If $p=3$, then $\th'$ and $T$ generate a 
Sylow 3-subgroup of $G$, say $P$. We claim that in this case the eigenspace
$$
L_\th=\{f\in L\mid\th f=(\det\th)f\}
$$
is a 1-dimensional subspace of $L$ spanned by a reducible polynomial, and 
therefore $t^S\notin L_\th$. Since Sylow 3-subgroups of $G$ are conjugate to 
each other and since $L$ is stable under all linear operators representing 
elements of $G$, it suffices to check the claim assuming that $P$ is one 
particular Sylow 3-subgroup.  Moreover, since all elements of $T$ are represented 
by linear operators acting on the whole $L$ as the multiplications by their 
determinants, it suffices to consider any particular element $\th'\in P$ not 
lying in $T$. So we take $\th'$ represented by the linear operator
$$
x_1\mapsto\ep x_1,\quad x_2\mapsto x_2,\quad x_3\mapsto x_3
$$
where $\ep$ is a primitive cube root of 1. Then $\det\th=\ep$ and 
$L_\th=\<x_1x_2x_3\>$. Indeed, $\th$ has two eigenvalues on $L$ since 
$x_1^3+x_2^3+x_3^3$ is fixed by $\th$. Thus we get a contradiction both for 
$p=2$ and for $p=3$.
\endproof

\proclaim
Lemma 4.2.
All elements of order $2$ in $G$ are conjugate to each other, and the same 
holds for the elements of order $4$. The nonidentity elements of $T$ are 
conjugate to each other in the group $G$.
\endproclaim

\Proof.
We view elements $g\in\SL_2(\bbF_3)$ as linear transformations of the vector 
space $\bbF_3^2$. If $g$ has order 2, then $g^2=\Id$. Such an operator is 
diagonalizable with eigenvalues equal to $-1$ or 1. Since $g\ne\Id$ and $\det 
g=1$, we must have $g=-\Id$. In other words, there is only one element of 
order 2 in $\SL_2(\bbF_3)$. Suppose that $g$ has order 4. Then $g^2=-\Id$. 
Taking any nonzero vector $w_1\in\bbF_3^2$ and $w_2=gw_1$, we get $gw_2=-w_1$. 
This means that $g$ has matrix
$$
\pmatrix{0&-1\cr 1&0}
$$
in a suitable basis of $\bbF_3^2$. Any two such operators are conjugate in 
$\GL_2(\bbF_3)$, but in fact they are conjugate in $\SL_2(\bbF_3)$. This 
follows from the fact that the centralizer of $g$ in $\GL_2(\bbF_3)$ contains 
a linear operator with determinant $-1$, e.g., the operator with matrix
$$
\pmatrix{1&-1\cr 1&1}
$$
in the previously considered basis of $\bbF_3^2$. Thus the elements of order 
4 form one conjugacy class in $\SL_2(\bbF_3)$.

Let $\pi:G\to\SL_2(\bbF_3)$ be the homomorphism with kernel $T$. If $g\in G$ 
has order 2 or 4, then $\pi(g)$ has the same order. Since no nonzero vector in 
$\bbF_3^2$ is fixed by $\pi(g)$, we have $gtg^{-1}\ne t$, and so $gt\ne tg$, 
for all $t\in T$. But then the map $t\mapsto t^{-1}gt$ gives a bijection of 
$T$ onto $gT$. Hence the conjugacy class of $g$ contains the coset $gT$. It 
follows that $g$ is conjugate to some element $h\in G$ if and only if $\pi(g)$ 
is conjugate to $\pi(h)$. Thus in $G$ there is one conjugacy class containing 
all elements of order 2 and one conjugacy class containing all elements of 
order 4.

The claim concerning conjugacy of elements of $T$ follows from the fact that 
$\SL_2(\bbF_3)$ acts transitively on nonzero vectors of $\bbF_3^2$.
\endproof

\proclaim
Corollary 4.3.
Each element $\th'\in G_A$ has order $\le4,$ and if the order of $\th'$ is 
$3,$ then $\th'\in T$.
\endproclaim

\Proof.
If $\th'\notin T$, then it follows from Proposition 4.1 that the coset 
$\th'T\in G/T$ has order 2 or 4. Since $|T|=9$ is prime to 2, such a coset 
contains an element of order 2 or 4. But we have seen in the proof of Lemma 
4.2 that all elements in such a coset belong to one conjugacy class, and so 
they have the same order.
\endproof

\proclaim
Lemma 4.4.
The Hessian group $G$ acts on the projective plane $\bbP^2$ in such a way that 
any linear operator $\tau\in\GL(V)$ representing an element $\tau'\in G$ 
extends  to an isomorphism $\Skl_3(a,b,c)\to\Skl_3(\at,\bt,\ct)$ whenever 
$\tau'(a:b:c)=(\at:\bt:\ct)$. If $R$ is a Hecke symmetry on $V$ such that 
$\bbS(V,R)=\Skl_3(a,b,c),$ then 
$$
R_\tau=(\tau\ot\tau)\circ R\circ(\tau^{-1}\ot\tau^{-1})
$$
is a Hecke symmetry satisfying $\bbS(V,R_\tau)=\Skl_3(\at,\bt,\ct),$ and if 
$\th$ and $\th_\tau$ are the linear operators given by Proposition 3.1 for $R$ 
and for $R_\tau,$ then $\th_\tau=\tau\th\tau^{-1}$.
\endproclaim

\Proof.
Let $W\sbs\bbT_3(V)$ be the subspace with a basis formed by the 3 tensors
$$
w_1=\sum x_{i-1}x_ix_{i+1},\qquad w_2=\sum x_{i+1}x_ix_{i-1},\qquad 
w_3=\sum x_i^3.
$$
We have $W=W^+\oplus W^-$ where $W^+$ is the intersection of $W$ with the 
subspace of symmetric tensors, and $W^-$ is the 1-dimensional subspace of 
alternating tensors in $\bbT_3(V)$. The subspace $W^+$ is mapped 
isomorphically onto the 2-dimensional subspace $L\sbs\bbS_3(V)$ spanned by the 
polynomials $x_1x_2x_3$ and $x_1^3+x_2^3+x_3^3$. Since $L$ is stable under the 
action of $\tau$ on $\bbS_3(V)$, while $W^-$ is a $\GL(V)$-submodule of 
$\bbT_3(V)$, we get $\tau^{\ot3}(W)=W$. Hence there is a homomorphism 
$G\to\PGL(W)$ which defines an action of $G$ on $\bbP^2=\bbP(W)$.

Now $a,b,c$ are the coordinates of the tensor $t=aw_1+bw_2+cw_3$ associated 
with the algebra $\Skl_3(a,b,c)$ in the chosen basis of $W$. We have 
$\tau'(a:b:c)=(\at:\bt:\ct)$ where $\at,\bt,\ct$ are the coordinates of the 
tensor $\tau^{\ot3}(t)\in W$. Since the space $J_2$ of defining relations of 
$\Skl_3(a,b,c)$ is the smallest subspace of $\bbT_2(V)$ with the property that 
$t\in J_2V$, the operator $\tau^{\ot2}$ maps it onto a similar subspace 
defined for $\at,\bt,\ct$. This means that $\tau$ extends to an isomorphism 
$\Skl_3(a,b,c)\to\Skl_3(\at,\bt,\ct)$.

Clearly $R_\tau$ satisfies the same Hecke relation and the braid relation 
satisfied by $R$, and $\tau$ extends to an isomorphism $\bbS(V,R)\to\bbS(V,R_\tau)$. 
Hence the defining relations of the algebra $\bbS(V,R_\tau)$ are those of 
$\Skl_3(\at,\bt,\ct)$. Finally, the equality $\th_\tau=\tau\th\tau^{-1}$ 
is deduced from the fact that the operator $\tau^{\ot4}$ intertwines the two 
representations of $\calH_4$ in $\bbT_4(V)$ arising from $R$ and from $R_\tau$.
\endproof

\Remark.
If $\tau\in\GL(V)$ is such that $\tau\ot\tau$ commutes with $R$, then the 
$\tau$-twist of $R$ is the Hecke symmetry 
$(\tau\ot\Id_V)\circ R\circ(\tau^{-1}\ot\Id_V)$. In \cite{Ew-O94} Hecke 
symmetries were determined up to a twist. However, we will not use such twists.
\endremark

\section
5. Algebras of type A are not associated with Hecke symmetries

The final result announced in the introduction is probably true over an 
arbitrary base field. However, it will be proved under the assumption that 
$\chr\bbk\ne 2,3$. When $\chr\bbk=2$ or $3$, there are differences in the 
description of automorphisms of the Artin-Schelter regular algebras of type 
$A$. Moreover, over a field of characteristic 3 the elliptic curve with the 
$j$-invariant 0 does not admit a Hessian normal form, and therefore 
Artin-Schelter regular algebras with such a point scheme are not realized in 
the family of algebras $\Skl_3(a,b,c)$.

\proclaim
Theorem 5.1.
Let $R$ be a Hecke symmetry on a $3$-dimensional vector space $V$ over a field\/ 
$\bbk$ of characteristic $\ne2,3$. Then $\bbS(V,R)$ cannot be any Artin-Schelter 
regular algebra of type $A$.
\endproclaim

\Proof.
Suppose that $\bbS(V,R)=\Skl_3(a,b,c)$ in some fixed basis $x_1,x_2,x_3$ of $V$, 
and assume that this algebra is Artin-Schelter regular of type $A$. Extending 
the base field $\bbk$ we may assume it to be algebraically closed.

Recall that the space of defining relations of the $R$-symmetric algebra 
$\bbS(V,R)$ is $\Up^{(2)}=\Im(R-q\cdot\Id)$. Hence
$$
\Up^{(2)}=\<t_1,t_2,t_3\>,\qquad\Up^{(3)}=
(\Up^{(2)}\ot V)\cap(V\ot\Up^{(2)})=\<t\>
$$
where $t_1,t_2,t_3$ and $t$ are as in section 4.

The operator $\th$ introduced in Proposition 3.1 extends to an automorphism of 
$\Skl_3(a,b,c)$, and by Proposition 4.1 its image $\th'$ in $\PGL(V^*)$ lies 
in the Hessian group $G$. By Lemma 4.4 there exists then a similar Hecke 
symmetry with a changed triple of parameters $a,b,c$ and with $\th'$ changed 
to any element in the same conjugacy class of $G$. Therefore we may assume 
$\th'$ to be any representative of this conjugacy class. By Lemma 4.2 and 
Corollary 4.3 there are at most 4 conjugacy classes of $G$ that can contain 
$\th'$. Accordingly, there are 4 cases to consider.

In each case we first use Proposition 3.9 to determine possible values of $q$. 
It will be seen that $q\ne-1$, i.e., $[2]_q\ne0$. Then we will use the relations 
between $t,f,\ph,\psi,\th$ established in Theorem 3.8 and Proposition 3.13. 
Note that $\psi=\Id$, and therefore $\ph=q^4\mskip1mu\th^{-2}$. The operator 
$\th$ is expressed as follows:
$$
\th(x_j)=q\,\sum_if(x_jt_i)\,x_i.
$$
Any Hecke symmetry with the prescribed algebra $\bbS(V,R)$ is completely 
determined by the function $f:V^{\ot3}\to\bbk$ which is not yet known. 
However, for each possible $\th$ we find the values $f(x_jt_i)$ from the 
previous expression. Combining them with the twisted cyclicity condition
$$
f(v_1v_2v_3)=f\bigl(\ph(v_3)v_1v_2\bigr),\qquad v_1,v_2,v_3\in V,
$$
we get a set of linear equations for the values of $f$ at the degree 3 
monomials in $x_1,x_2,x_3$. The functions obtained by solving these equations 
may depend on several additional parameters. However, fulfillment of all 
linear equations arising from Proposition 3.13 does not guarantee the 
existence of a Hecke symmetry, and we have to look further.

By Proposition 3.12 the generating space $I_2$ of the ideal $I\sbs\bbT(V)$ 
defining the algebra $\La(V,R)$ is determined by $f$ as
$$
I_2=\{u\in\bbT_2(V)\mid f(uV)=0\}=\{u\in\bbT_2(V)\mid f(Vu)=0\}.
$$
The projection $P$ of $V^{\ot2}$ onto $\Up^{(2)}$ with respect to the direct 
sum decomposition $V^{\ot2}=I_2\oplus\Up^{(2)}$ is expressed as
$$
P(w)=\sum f(\xwt_iw)t_i,\qquad w\in V^{\ot2},
$$
where $\xwt_1,\xwt_2,\xwt_3$ is the basis of $V$ dual to the basis 
$t_1,t_2,t_3$ of $\Up^{(2)}$ with respect to the pairing 
$V\times\Up^{(2)}\to\bbk$ induced by $f$, i.e., $f(\xwt_it_j)$ is 1 for $i=j$ 
and 0 otherwise. Since $I_2$ and $\Up^{(2)}$ are the eigenspaces of $R$ 
corresponding to the eigenvalues $q$ and $-1$, we get
$$
R=-P+q\,(\Id-P)=q\Id{}-(1+q)P.
$$
The braid equation for $R$ is equivalent to the following equation for $P$:
$$
\openup1\jot
\eqalign{
&(P\ot\Id_V)(\Id_V\ot\,P)(P\ot\Id_V)-{q\over(1+q)^2}(P\ot\Id_V)\cr
&\hskip3cm{}=(\Id_V\ot\,P)(P\ot\Id_V)(\Id_V\ot\,P)-{q\over(1+q)^2}(\Id_V\ot\,P).
}
$$
Note that $P\ot\Id$ and $\Id\ot\,P$ are projections (idempotent linear 
operators) mapping $V^{\ot3}$ onto $\Up^{(2)}\ot V=\Up^{(2,1)}$ and 
$V\ot\Up^{(2)}=\Up^{(1,2)}$, respectively. If the displayed equality is 
true, then the operators in the left and right hand sides of the formula must 
have images in $\Up^{(2,1)}\cap\Up^{(1,2)}=\Up^{(3)}$. Thus
$$
(\Id\ot\,P)(P\ot\Id)u-q(1+q)^{-2}u\in\Up^{(3)}\qquad\hbox{for all 
$u\in\Up^{(1,2)}$}.
$$
This can be rephrased by saying that $\Id\ot\,P$ and $P\ot\Id$ induce linear 
maps between the vector spaces $\Up^{(2,1)}/\Up^{(3)}$ and 
$\Up^{(1,2)}/\Up^{(3)}$ which are inverse to each other up to a scalar 
multiple. This condition leads to a set of quadratic equations on the 
parameters involved in the construction of $f$. It will be shown that in each 
case no solution is possible, and so the braid equation for $R$ cannot be 
satisfied.

\medskip
{\bf Case 1:} $\th'$ is the identity element. In this case $\th=\la\Id$ for 
some $\la\in\bbk\setm\{0\}$. Then $\tr\th=3\la$, and Proposition 3.9 yields 
$3\la=q(1+q+q^2)$. By Corollary 3.11 $\la=q^2$. Hence $q(1-2q+q^2)=0$, and we 
find $q=1$. It follows that $\th=\Id$, and $\ph=\Id$ as well.

The expression of $\th$ in terms of the function $f$ gives $f(x_jt_i)=0$ for 
$i\ne j$, and $f(x_it_i)=1$ for all $i$. On the other hand, the cyclicity 
property $f(v_1v_2v_3)=f(v_3v_1v_2)$ allows us to write
$$
\openup1\jot
\eqalign{
f(x_{i-1}t_i)&=af(x_{i-1}x_{i+1}x_{i-1})+bf(x_{i-1}^2x_{i+1})+cf(x_{i-1}x_i^2)\cr 
&=(a+b)f(x_{i+1}x_{i-1}^2)+cf(x_{i-1}x_i^2).
}
$$
The equalities $f(x_{i-1}t_i)=0$ give a system of linear equations for the 
values of $f$ at 3 monomials:
$$
\openup1\jot
\eqalign{
(a+b)f(x_2x_3^2)+cf(x_3x_1^2)=0,\cr
(a+b)f(x_3x_1^2)+cf(x_1x_2^2)=0,\cr
(a+b)f(x_1x_2^2)+cf(x_2x_3^2)=0.
}
$$
Its matrix of coefficients has determinant
$$
\left|\matrix{a+b & c & 0\cr
0 & a+b & c \cr
c & 0 & a+b}\right|=(a+b)^3+c^3\ne0.
$$
Hence the system admits only the zero solution, i.e., $f(x_{i-1}x_i^2)=0$ for 
all $i$. Similarly, the equalities $f(x_{i+1}t_i)=0$ force $f(x_{i+1}x_i^2)=0$ 
for all $i$.

By cyclicity of $f$ we have $\,f(x_{i-1}x_ix_{i+1})=a'\,$ and 
$\,f(x_{i+1}x_ix_{i-1})=b'\,$ for some $a',b'\in\bbk$ which do not depend on 
$i$. The equalities $f(x_it_i)=1$ are now written as
$$
aa'+bb'+cf(x_i^3)=1.
$$
It follows that $f(x_i^3)=c'$ for all $i$ where $c'\in\bbk$ is such that 
$aa'+bb'+cc'=1$.

The projection $P$ is given by the formula $\,P(w)=\sum f(x_iw)t_i\,$ for 
$w\in V^{\ot2}$, i.e.,
$$
P(x_{i+1}x_{i-1})=a't_i,\qquad P(x_{i-1}x_{i+1})=b't_i,\qquad P(x_i^2)=c't_i.
$$
The projections $\Id\ot\,P$ and $P\ot\Id$ induce linear maps between 3 pairs 
of subspaces of $\Up^{(2,1)}$ and $\Up^{(1,2)}$. One pair is formed by the 
subspaces $\<t_3x_1,t_1x_2,t_2x_3\>\sbs\Up^{(2,1)}$ and 
$\<x_1t_3,x_2t_1,x_3t_2\>\sbs\Up^{(1,2)}$, and the matrices of linear maps in 
the respective bases are
$$
\pmatrix{ab' & ca' & bc' \cr bc' & ab' & ca' \cr ca' & bc' & ab' },\qquad
\pmatrix{ba' & cb' & ac' \cr ac' & ba' & cb' \cr cb' & ac' & ba' }.
$$
Since $\<x_1t_3,x_2t_1,x_3t_2\>\cap\Up^{(3)}=0$, the composite map 
$(\Id\ot\,P)(P\ot\Id)$, when restricted to $\<x_1t_3,x_2t_1,x_3t_2\>$, must be 
a scalar multiplication. This leads to equations
$$
\openup1\jot
\eqalign{
bcb'c'+cac'a'+aba'b'&=\ka,\cr
bca'^2+cab'^2+abc'^2&=0,\cr
a^2b'c'+b^2c'a'+c^2a'b'&=0
}
$$
where $\ka=q(1+q)^{-2}=1/4$ (the precise value of $\ka$ turns out to be 
irrelevant). The second pair of subspaces $\<t_2x_1,t_3x_2,t_1x_3\>$ and 
$\<x_1t_2,x_2t_3,x_3t_1\>$ produces the same result. Finally, the restrictions 
of $\Id\ot\,P$ and $P\ot\Id$, respectively, to the subspaces 
$\<t_1x_1,t_2x_2,t_3x_3\>$ and $\<x_1t_1,x_2t_2,x_3t_3\>$ give linear maps 
with the matrices
$$
\pmatrix{cc' & bb' & aa' \cr aa' & cc' & bb' \cr bb' & aa' & cc' },\qquad
\pmatrix{cc' & aa' & bb' \cr bb' & cc' & aa' \cr aa' & bb' & cc' }.
$$
Since $\,(\Id\ot\,P)(P\ot\Id)(x_it_i)-\ka\,x_it_i\,$ must be a scalar multiple 
of $t$, we get
$$
a^2a'^2+b^2b'^2+c^2c'^2-\ka=bcb'c'+cac'a'+aba'b'.
$$
Combining the last equation with the previous three, we arrive at a system of 
three homogeneous quadratic equations
$$
\openup1\jot
\eqalign{
bc\,a'^2+ca\,b'^2+ab\,c'^2&=0,\cr
a^2b'c'+b^2c'a'+c^2a'b'&=0,\cr
a^2a'^2+b^2b'^2+c^2c'^2&=2bc\,b'c'+2ca\,c'a'+2ab\,a'b'
}
$$
which we aim to solve in $a',b',c'$. Note that these are the same equations 
as in the two papers of Ohn \cite{Ohn99, (10.5)} and \cite{Ohn05, (5)}. 
A nonzero solution exists if and only if the resultant of the system vanishes. 
The resultant of the three polynomials
$$
\openup1\jot
\eqalign{
F_1&{}=bc\,X^2+ca\,Y^2+ab\,Z^2,\cr
F_2&{}=a^2YZ+b^2ZX+c^2XY,\cr
F_3&{}=a^2X^2+b^2Y^2+c^2Z^2-2bc\,YZ-2ca\,ZX-2ab\,XY
}
$$
can be computed as a certain determinant of order 6 by the formula of Sylvester 
(see \cite{Gel-KZ, Ch. 3, section 4D}). We need three other quadratic 
polynomials $D_1$, $D_2$, $D_3$ determined uniquely only modulo the linear 
span of $F_1$, $F_2$, $F_3$. One has
$$
D_1=\det(l_{ij})_{1\le i,j\le3}
$$
where $l_{1j}\in\bbk$ and $l_{2j},\,l_{3j}$ are linear forms such that
$F_j=l_{1j}X^2+l_{2j}Y+l_{3j}Z$. In the case considered here we may take
$$
\openup2\jot
\displaylines{
D_1=\left|\matrix{
bc & 0 & a^2 \cr
acY & c^2X+a^2Z & b^2Y-2abX-2bcZ \cr
abZ & b^2X & c^2Z-2acX \cr
}\right|\cr
=2abc(b^3-c^3)X^2+a^2b(c^3-a^3)Z^2+b^2c(a^3-b^3)XY+
bc^2(2b^3+c^3-3a^3)XZ.
}
$$
The polynomials $D_2$ and $D_3$ are defined similarly using expressions of 
$F_j$, respectively, as $l_{1j}X+l_{2j}Y^2+l_{3j}Z$ and as 
$l_{1j}X+l_{2j}Y+l_{3j}Z^2$. For $D_2$, $D_3$ we take
$$
\openup1\jot
\eqalign{ 
&2abc(c^3-a^3)Y^2+b^2c(a^3-b^3)X^2+c^2a(b^3-c^3)YZ+
ca^2(2c^3+a^3-3b^3)YX,\cr
&2abc(a^3-b^3)Z^2+c^2a(b^3-c^3)Y^2+a^2b(c^3-a^3)ZX+
ab^2(2a^3+b^3-3c^3)ZY.
}
$$
By Sylvester's formula the resultant $\Res(F_1,F_2,F_3)$ is equal to the 
determinant
$$
\left|\matrix{
bc&0&a^2&2abc(b^3-c^3)&b^2c(a^3-b^3)&0\cr
ac&0&b^2&0&2abc(c^3-a^3)&ac^2(b^3-c^3)\cr
ab&0&c^2&a^2b(c^3-a^3)&0&2abc(a^3-b^3)\cr
0&a^2&-2bc&0&ac^2(b^3-c^3)&ab^2(2a^3+b^3-3c^3)\cr
0&b^2&-2ac&bc^2(2b^3+c^3-3a^3)&0&a^2b(c^3-a^3)\cr
0&c^2&-2ab&b^2c(a^3-b^3)&a^2c(2c^3+a^3-3b^3)&0\cr
}\right|
$$
whose columns are composed of the coefficients of $X^2$, $Y^2$, $Z^2$, $YZ$, 
$ZX$, $XY$ in $F_1$, $F_2$, $F_3$, $D_1$, $D_2$, $D_3$. A machine computation 
gives the following result:
$$
\displaylines{
a^2b^2c^2(a^{18}+b^{18}+c^{18})
+6\,a^2b^2c^2(a^{15}b^3+a^{15}c^3+a^3b^{15}+a^3c^{15}+b^{15}c^3+b^3c^{15})\cr
+15\,a^2b^2c^2(a^{12}b^6+a^{12}c^6+a^6b^{12}+a^6c^{12}+b^{12}c^6+b^6c^{12})\cr
+20\,a^2b^2c^2(a^9b^9+a^9c^9+b^9c^9)
-24\,a^5b^5c^5(a^9+b^9+c^9)\cr
-102\,a^5b^5c^5(a^6b^3+a^6c^3+a^3b^6+a^3c^6+b^6c^3+b^3c^6)
+495\,a^8b^8c^8.\cr
}
$$
Expressing this in contracted form, we conclude that
$$
\Res(F_1,F_2,F_3)=a^2b^2c^2\bigl((a^3+b^3+c^3)^3-27a^3b^3c^3\bigr)^2\ne0.
$$
Thus we do not get any solution.

\medskip
{\bf Case 2:} $\th'$ has order 3. We may assume that
$$
\th(x_i)=\la\ep^ix_i,\qquad i=1,2,3,
$$
where $\ep$ is a primitive cube root of 1 and $\la\in\bbk\setm\{0\}$. Since 
$\tr\th=0$, Proposition 3.9 shows that $1+q+q^2=0$, i.e., $q$ is a primitive 
cube root of 1 too.

Since $\th^{\ot3}(t)=\la^3t$, the identity $\th^{\ot3}(t)=q^6t$ yields 
$\la^3=q^6$. Hence $\la=q^2\om$ where $\om\in\{1,\ep,\ep^2\}$. Renumbering 
the generators $x_i$, we may assume that $\om=1$, and so $\la=q^2$. Then
$$
\ph(x_i)=q^4\th^{-2}(x_i)=\ep^{-2i}x_i=\ep^ix_i.
$$
The cyclicity of $f$ is expressed as $f(wx_i)=\ep^if(x_iw)$ for 
$w\in\bbT_2(V)$. In particular, $f(x_i^3)=\ep^if(x_i^3)$ and 
$f(x_ix_jx_k)=\ep^{i+j+k}f(x_ix_jx_k)$. It follows that
$$
f(x_1^3)=f(x_2^3)=0,\qquad f(x_ix_jx_k)=0\quad\hbox{whenever 
$i+j+k\not\equiv0\mod 3$}.
$$
Let\quad$f(x_1x_2x_3)=f(x_3x_1x_2)=qa'$,\quad 
$f(x_2x_1x_3)=f(x_3x_2x_1)=qb'$,\quad $f(x_3^3)=qc'$.
Then $f(x_2x_3x_1)=\ep qa'$ and $f(x_1x_3x_2)=\ep^2qb'$.

The formula expressing $\th$ in terms of $f$ shows that $f(x_jt_i)=0$ whenever 
$i\ne j$, and $f(x_it_i)=q\ep^i$ for each $i$. Therefore the projection $P$ 
is written as
$$
P(w)=q^{-1}\sum\ep^{-i}f(x_iw)t_i,\qquad w\in\bbT_2(V).
$$
So
$$
\vcenter{\halign{
\hfil$#$&${}=#$\hfil & \qquad\hfil$#$&${}=#$\hfil & \qquad
\hfil$#$&${}=#$\hfil \cr 
P(x_1^2) & 0, & P(x_2x_3) & \ep^2a't_1, & P(x_3x_2) & \ep b't_1, \cr
P(x_2^2) & 0, & P(x_3x_1) & \ep^2a't_2, & P(x_1x_3) & \ep b't_2, \cr
P(x_3^2) & c't_3, & P(x_1x_2) & a't_3, & P(x_2x_1) & b't_3.\cr
}}
$$
The linear maps between $\<t_3x_1,t_1x_2,t_2x_3\>$ and 
$\<x_1t_3,x_2t_1,x_3t_2\>$ induced by the projections $\Id\ot\,P$ and 
$P\ot\Id$ have matrices
$$
\pmatrix{ab' & ca' & bc' \cr 0 & \ep ab' & \ep^2ca' \cr \ep^2ca' & 0 & \ep ab'}
\qquad\hbox{and}\qquad
\pmatrix{ba' & cb' & ac' \cr 0 & \ep^2ba' & \ep cb' \cr \ep cb' & 0 & \ep^2ba'}.
$$
We compute
$$
(\Id\ot\,P)(P\ot\Id)(x_1t_3)=(aba'b'+\ep bcb'c')x_1t_3+c^2a'b'x_2t_1
+\ep^2(bca'^2+cab'^2)x_3t_2.
$$
Since $(\Id\ot\,P)(P\ot\Id)|_{\<x_1t_3,x_2t_1,x_3t_2\>}$ must be the 
multiplication by $q(1+q)^{-2}$, we deduce that $bca'^2+cab'^2=c^2a'b'=0$. 
Therefore $a'=b'=0$ since $abc\ne0$. But then
$$
(\Id\ot\,P)(P\ot\Id)(x_1t_3)=0\ne q(1+q)^{-2}x_1t_3\,,
$$
a contradiction.
%
%
%

\medskip
{\bf Case 3:} $a=b$ and $\th'$ has order 2. We may assume that
$$
\th(x_1)=\la x_2,\qquad\th(x_2)=\la x_1,\qquad\th(x_3)=\la x_3
$$
for some $\la\in\bbk\setm\{0\}$. Then $\,t_1\mapsto\la^2t_2$, $\,t_2\mapsto\la^2t_1$, 
$\,t_3\mapsto\la^2t_3$ under the action of $\th^{\ot2}$. We see that 
$\,\tr\th=\la\,$ and $\,\tr\th^{\ot2}\mid_{\Up^{(2)}}=\la^2\,$, whence
$$ 
\la=q(1+q+q^2),\qquad\la^2=q^3(1+q+q^2)
$$
by Proposition 3.9. It follows that $\,\la=q^2=-1\,$ and $\,\ph=q^4\,\th^{-2}=\Id$.

The expression of $\th$ in terms of the function $f$ yields
$$
f(x_1t_2)=f(x_2t_1)=f(x_3t_3)=q^{-1}\la=q
$$
and $f(x_jt_i)=0$ whenever $(i,j)\ne(1,2),\,(2,1),\,(3,3)$. Using the cyclicity 
of $f$ as in Case 1, we get
$$
\openup1\jot
\eqalign{
2af(x_3x_1^2)+cf(x_1x_2^2)&{}=f(x_1t_2)=q,\cr
2af(x_1x_2^2)+cf(x_2x_3^2)&{}=f(x_2t_3)=0,\cr
2af(x_2x_3^2)+cf(x_3x_1^2)&{}=f(x_3t_1)=0
}
$$
(note that $2a=a+b$). This is a system of linear equations in 3 unknown values 
of $f$. Its matrix of coefficients has determinant
$$
d=8a^3+c^3=(a+b)^3+c^3\ne0.
$$
Hence the system admits a unique solution
$$
f(x_2x_3^2)=-2qacd^{-1},\qquad f(x_3x_1^2)=4qa^2d^{-1},\qquad
f(x_1x_2^2)=qc^2d^{-1}.
$$
In a similar way the values $f(x_2t_1)=q$, $f(x_1t_3)=f(x_3t_2)=0$ lead to a 
system of equations yielding
$$
f(x_1x_3^2)=-2qacd^{-1},\qquad f(x_3x_2^2)=4qa^2d^{-1},\qquad
f(x_2x_1^2)=qc^2d^{-1}.
$$
Next, setting $f(x_1x_2x_3)=qa'$ and $f(x_2x_1x_3)=qb'$, we get
$$
\openup1\jot
\eqalign{
qa(a'+b')+cf(x_1^3)&{}=f(x_1t_1)=0,\cr
qa(a'+b')+cf(x_2^3)&{}=f(x_2t_2)=0,\cr
qa(a'+b')+cf(x_3^3)&{}=f(x_3t_3)=q\,.
}
$$
It follows that
$$
f(x_1^3)=f(x_2^3)=qc',\qquad f(x_3^3)=qc''
$$
for $c',c''\in\bbk$ such that\quad $a(a'+b')+cc'=0$\quad and\quad $c''=c'+c^{-1}$.

The projection $P$ with $\Im P=\Up^{(2)}$ and $\Ker P=I_2$ is given by
$$
P(w)=q^{-1}\bigl(f(x_2w)t_1+f(x_1w)t_2+f(x_3w)t_3\bigr).
$$
So
$$
\eqalign{
P(x_1^2)&{}=d^{-1}(c^2t_1+dc't_2+4a^2t_3),\cr
P(x_2^2)&{}=d^{-1}(dc't_1+c^2t_2+4a^2t_3),\cr
P(x_3^2)&{}=d^{-1}(-2act_1-2act_2+dc''t_3),
}
$$
$$
\openup1\jot
\vcenter{\halign{\hfil$#$&${}=#$\hfil & \qquad\hfil$#$&${}=#$\hfil \cr 
P(x_2x_3) & d^{-1}(4a^2t_1+da't_2-2act_3),
& P(x_3x_2) & d^{-1}(4a^2t_1+db't_2-2act_3), \cr
P(x_3x_1) & d^{-1}(da't_1+4a^2t_2-2act_3), 
& P(x_1x_3) & d^{-1}(db't_1+4a^2t_2-2act_3), \cr
P(x_1x_2) & d^{-1}(c^2t_1+c^2t_2+da't_3), 
& P(x_2x_1) & d^{-1}(c^2t_1+c^2t_2+db't_3).\cr
}}
$$
The projection $\Id\ot\,P$ gives by restriction a linear map 
$\Up^{(2,1)}\to\Up^{(1,2)}$. Its matrix with respect to the bases 
$t_1x_1,t_2x_1,t_3x_1,\ldots,t_3x_3$ and 
$x_1t_1,x_1t_2,x_1t_3,\ldots,x_3t_3$ of those two 
spaces is $d^{-1}M$ where
$$
\hfuzz=10pt
M\!=\!\pmatrix{
c^3   & daa'  & ac^2   & c^3   & 4a^3   & dac'   & dcb'   & -2a^2c & 4a^3   \cr
dcc'  & 4a^3  & ac^2   & c^3   & dab'   & ac^2   & 4a^2c  & -2a^2c & daa'   \cr
4a^2c & -2a^2c& dab'   & dca'  & -2a^2c & 4a^3   & -2ac^2 & dac''  & -2a^2c \cr
daa'  & c^3   & ac^2   & 4a^3  & dcc'   & ac^2   & -2a^2c & 4a^2c  & dab'   \cr
4a^3  & c^3   & dac'   & dab'  & c^3    & ac^2   & -2a^2c & dca'   & 4a^3   \cr
-2a^2c& dcb'  & 4a^3   & -2a^2c& 4a^2c  & daa'   & dac''  & -2ac^2 & -2a^2c \cr
ac^2  & ac^2  & dca'   & dac'  & ac^2   & 4a^2c  & 4a^3   & dab'   & -2ac^2 \cr
ac^2  & dac'  & 4a^2c  & ac^2  & ac^2   & dcb'   & daa'   & 4a^3   & -2ac^2 \cr
dab'  & 4a^3  & -2ac^2 & 4a^3  & daa'   & -2ac^2 & -2a^2c & -2a^2c & dcc''
}
$$
The other projection $P\ot\Id$ gives a map in the opposite direction with the 
matrix $d^{-1}N$ where $N$ is obtained from $M$ by interchanging $a'$ and 
$b'$.
%
%
%

This interrelation between $M$ and $N$ is explained by the fact that 
$d^{-1}N$ coincides with the matrix of the linear map 
$\Up^{(2,1)}\to\Up^{(1,2)}$ obtained by restriction of the linear operator 
$\Id_V\ot\,\tau P\tau$ where $\tau:V^{\ot2}\to V^{\ot2}$ is the flip of 
tensorands. Note that all elements of $\Up^{(2)}$ are fixed by $\tau$. The 
linear operator $\tau P\tau$ is the projection of $V^{\ot2}$ onto $\Up^{(2)}$ 
corresponding to the Hecke symmetry $R\op=\tau R\tau$, and the latter is 
associated with the linear function $f\op$ on $V^{\ot3}$ whose values on the 
monomials in $x_1,x_2,x_3$ are given by the same formulas as the values of 
$f$, but with $a'$ and $b'$ interchanged. This observation explains also that 
whenever some relation between the parameters $a,c,a',b',c'$ must hold to 
satisfy the Hecke symmetry conditions, the companion relation with $a'$ and 
$b'$ interchanged must hold too.

Since $(\Id\ot\,P)(P\ot\Id)u\equiv q(1+q)^{-2}u\mod{\Up^{(3)}}$ for all 
$u\in\Up^{(1,2)}$, the product $MN$ must have zero entries everywhere 
except the 1st, 5th, 9th rows and the principal diagonal. We will only need a 
few values in the 1st column of the matrix $MN$. Computing its entries in 
positions (2,1), (4,1), (7,1) and (8,1), we get
$$
\displaylines{
\hfill d(cc'+ab')(c^3+4a^3)+4a^3c^3+(daa')^2=0,\hfill\llap{(1)}\cr
\hfill d(aa'+cc')c^3+4a^3c^3+4a^3d(ab'+cc')+(dab')(daa')=0,\hfill\llap{(2)}\cr
\hfill ac^5+ac^2d(cc'+2aa'+ab')+d^2a^2c'b'=0,\hfill\llap{(3)}\cr
\hfill ac^5+d^2acc'^2+24a^4c^2+ac^2(-dab'-daa')=0.\hfill\llap{(4)}
}
$$
Making substitutions $cc'+ab'=-aa'$ and $aa'+cc'=-ab'$ in the first two 
equalities, we rewrite them as
$$
(daa'-c^3)(daa'-4a^3)=0,\qquad (daa'-c^3)(dab'-4a^3)=0.
$$
As we have noted, the equalities remain true with $a'$ and $b'$ interchanged. 
Thus
$$
(dab'-c^3)(dab'-4a^3)=0,\qquad (dab'-c^3)(daa'-4a^3)=0.
$$
It follows that either $daa'=dab'=c^3$ or $daa'=dab'=4a^3$. In any case 
$a'=b'$ since $d\ne0$ and $a\ne0$. Multiplying (3) by $a^{-1}c$ and 
substituting $b'=a'$, $cc'=-2aa'$, we arrive at
$$
(c^3-daa')(c^3+2daa')=0.
$$
If $daa'\ne c^3$, then $daa'=4a^3$, but in this case $c^3+2daa'=c^3+8a^3=d\ne0$, 
and we get a contradiction. Hence $daa'=dab'=c^3$ and $dcc'=-2c^3$. But then 
(4) reduces to $3ac^5+24a^4c^2=0$, i.e., $3ac^2d=0$, again in contradiction 
with $ac\ne0$ and $d\ne0$.

\medskip
{\bf Case 4:} $a=b$ and $\th'$ has order 4. We may assume that
$$
\th(x_j)=\la\sum_{i=1}^3\ep^{ij}x_i,\qquad j=1,2,3,
$$
where $\la\in\bbk\setm\{0\}$ and $\ep$ is a primitive cube root of 1. Since 
$a=b$, the tensor $t$ is symmetric, and its image in the symmetric algebra of 
$V$ is the polynomial
$$
t^S=c(x_1^3+x_2^3+x_3^3)+6ax_1x_2x_3\in\bbS(V).
$$
Put $\ka=ac^{-1}$, so that $a=\ka c$. The condition that $\th$ is the degree 1 
component of an automorphism of $A$ means precisely that $\th$ sends $t^S$ to 
its scalar multiple. Since
$$
\th t^S=3(c+2a)\la^3(x_1^3+x_2^3+x_3^3)+18(c-a)\la^3x_1x_2x_3,
$$
we deduce that $\ka=(1-\ka)/(1+2\ka)$, i.e., $2\ka^2+2\ka=1$ (the curve defined 
by the equation $t^S=0$ has then the $j$-invariant 1728). In this case 
$\th t^S=(3+6\ka)\la^3t^S\!$, and therefore $\,\th^{\ot3}(t)=(3+6\ka)\la^3t\,$. 
The identity $\,\th^{\ot3}(t)=q^6t\,$ entails
$$
(3+6\ka)\la^3=q^6.
$$
Making use of relations $\ep^3=1$ and $1+\ep+\ep^2=0$, we get
$$
\eqalign{
\th^{\ot2}(t_j)&{}=a\,\th(x_{j+1})\th(x_{j-1})+a\,\th(x_{j-1})\th(x_{j+1})+c\,\th(x_j)^2\cr
&{}=\la^2c\sum_{m=1}^3\sum_{n=1}^3
\bigl(\ka\,\ep^{m(j+1)+n(j-1)}+\ka\,\ep^{m(j-1)+n(j+1)}+\ep^{(m+n)j}\bigr)x_mx_n\cr
&{}=\la^2c\sum_{i=1}^3\ep^{2ij}
\bigl((1-\ka)x_{i+1}x_{i-1}+(1-\ka)x_{i-1}x_{i+1}+(1+2\ka)x_i^2\bigr)\cr
&{}=\la^2(1+2\ka)\sum_{i=1}^3\ep^{2ij}t_i\,.
}
$$
Computation of traces gives $\,\tr\th=\la(1+2\ep)=\la(\ep-\ep^2)\,$ and
$$
\tr\th^{\ot2}|_{\Up^{(2)}}=\la^2(1+2\ka)(1+2\ep^2)=\la^2(1+2\ka)(\ep^2-\ep).
$$
By Proposition 3.9
$$
\la(\ep-\ep^2)=q(1+q+q^2),\qquad\la^2(1+2\ka)(\ep^2-\ep)=q^3(1+q+q^2).
$$
Since $\la\ne0$ and $\ep\ne\ep^2$, comparison of the last two equalities 
yields $(1+2\ka)\la=-q^2$. Since $(1+2\ka)^2=3$, it follows that
$$
(3+6\ka)\la^3=\bigl((1+2\ka)\la\bigr)^3=-q^6,
$$
which contradicts a relation found earlier.
\endproof

\references
\nextref Ab-An99
\auth{A.,Abella;N.,Andruskiewitsch}
\paper{Compact quantum groups arising from the FRT-construction}
\journal{Bol. Acad. Nac. Cienc. C\'ordoba}
\Vol{63}
\Year{1999}
\Pages{15-44}

\nextref And17
\auth{N.,Andruskiewitsch}
\paper{An introduction to Nichols algebras}
\InBook{Quantization, Geometry and Noncommutative Structures in Mathematics and Physics}
\publisher{Springer}
\Year{2017}
\Pages{135-195}

\nextref Ar-Sch87
\auth{M.,Artin;W.F.,Schelter}
\paper{Graded algebras of global dimension $3$}
\journal{Adv. Math.}
\Vol{66}
\Year{1987}
\Pages{171-216}

\nextref Ar-TV90
\auth{M.,Artin;J.,Tate;M.,Van den Bergh}
\paper{Some algebras associated to automorphisms of elliptic curves}
\InBook{The Grothendieck Festschrift, Volume I}
\publisher{Birkh\"auser}
\Year{1990}
\Pages{33-85}

\nextref Ar-TV91
\auth{M.,Artin;J.,Tate;M.,Van den Bergh}
\paper{Modules over regular algebras of dimension $3$}
\journal{Invent. Math.}
\Vol{106}
\Year{1991}
\Pages{335-388}

\nextref Bj-Br
\auth{A.,Bj\"orner;F.,Brenti}
\book{Combinatorics of Coxeter Groups}
\publisher{Springer}
\Year{2005}

\nextref Bond-P93
\auth{A.I.,Bondal;A.E.,Polishchuk}
\paper{Homological properties of associative algebras: The method of helices\inRus}
\journal{Izv. Ross. Akad. Nauk Ser. Mat.}
\Vol{57}
\nombre{2}
\Year{1993}
\Pages{3-50};
\etransl{Russian Acad. Sci. Izv. Math.}
\Vol{42}
\Year{1994}
\Pages{219-260}

\nextref Br-Kn
\auth{E.,Brieskorn;H.,Kn\"orrer}
\book{Plane Algebraic Curves}
\publisher{Birkh\"auser}
\Year{1986}

\nextref Ch-WW19
\auth{A.,Chirvasitu;C.,Walton;X.,Wang}
\paper{On quantum groups associated to a pair of preregular forms}
\journal{J.~Noncommut. Geom.}
\Vol{13}
\Year{2019}
\Pages{115-159}

\nextref Dub14
\auth{M.,Dubois-Violette}
\paper{Poincar\'e duality for Koszul algebras}
\InBook{Algebra, Geometry and Mathematical Physics}
\publisher{Springer}
\Year{2014}
\Pages{3-26}

\nextref Ew-O94
\auth{H.,Ewen;O.,Ogievetsky}
\paper{Classification of the $GL(3)$ quantum matrix groups}\hfill\break
\hbox{}\hfill arXiv:9412009.

\nextref Fi-MS97
\auth{D.,Fischman;S.,Montgomery;H.-J.,Schneider}
\paper{Frobenius extensions of subalgebras of Hopf algebras}
\journal{Trans. Amer. Math. Soc.}
\Vol{349}
\Year{1997}
\Pages{4857-4895}

\nextref Geck-P
\auth{M.,Geck;G.,Pfeiffer}
\book{Characters of Finite Coxeter Groups and Iwahori-Hecke Algebras}
\publisher{Clarendon Press}
\Year{2000}

\nextref Gel-KZ
\auth{I.M.,Gelfand;M.,Kapranov;A.,Zelevinsky}
\book{Discriminants, Resultants, and Multidimensional Determinants}
\publisher{Birkh\"auser}
\Year{1994}

\nextref Gur90
\auth{D.I.,Gurevich}
\paper{Algebraic aspects of the quantum Yang-Baxter equation\inRus}
\journal{Algebra i Analiz}
\Vol{2}
\nombre{4}
\Year{1990}
\Pages{119-148};
\etransl{Leningrad Math. J.}
\Vol{2}
\Year{1991}
\Pages{801-828}

\nextref Hai99
\auth{P.H.,Hai}
\paper{Poincar\'e series of quantum spaces associated to Hecke operators}
\journal{Acta Math. Vietnam}
\Vol{24}
\Year{1999}
\Pages{235-246}

\nextref Heck-Sch
\auth{I.,Heckenberger;H.-J.,Schneider}
\book{Hopf Algebras and Root Systems}
\publisher{Amer. Math. Soc.}
\Year{2020}

\nextref Lyu87
\auth{V.V.,Lyubashenko}
\paper{Vectorsymmetries}
\journal{Reports Dept. Math. Univ. Stockholm}
1987, No. 19.

\nextref Man
\auth{Yu.,Manin}
\book{Quantum Groups and Non-Commutative Geometry}
\publisher{CRM, Univ. Montr\'eal}
\Year{1988}

\nextref Ohn99
\auth{C.,Ohn}
\paper{Quantum $SL(3,\bbC)$'s with classical representation theory}
\journal{J.~Algebra}
\Vol{213}
\Year{1999}
\Pages{721-756}

\nextref Ohn05
\auth{C.,Ohn}
\paper{Quantum $SL(3,\bbC)$'s: the missing case}
\InBook{Hopf Algebras in Noncommutative Geometry and Physics}
\publisher{Marcel Dekker}
\Year{2005}
\Pages{245-255}

\nextref Pol-P
\auth{A.,Polishchuk;L.,Positselski}
\book{Quadratic Algebras}
\publisher{Amer. Math. Soc.}
\Year{2005}

\nextref Ros98
\auth{M.,Rosso}
\paper{Quantum groups and quantum shuffles}
\journal{Invent. Math.}
\Vol{133}
\Year{1998}
\Pages{399-416}

\nextref Skr20
\auth{S.,Skryabin}
\paper{On the graded algebras associated with Hecke symmetries}
\journal{J.~Noncommut. Geom.}
\Vol{14}
\Year{2020}
\Pages{937-986}

\nextref Sm96
\auth{S.P.,Smith}
\paper{Some finite dimensional algebras related to elliptic curves}
\InBook{Representation Theory of Algebras and Related Topics}
\publisher{Amer. Math. Soc.}
\Year{1996}
\Pages{315-348}

\endreferences
\bye